\documentclass[12pt]{article}
\usepackage{amssymb}
\newcommand{\R}{\mathbb{R}}

\newcommand{\Z}{\mathbb{Z}}
\renewcommand{\P}{\mathbb{P}}
\newcommand{\Q}{\mathbb{Q}}

\newcommand{\N}{\mathbb{N}}

\newcommand{\T}{\mathbb{T}}
\def\build#1_#2^#3{\mathrel{
\mathop{\kern 0pt#1}\limits_{#2}^{#3}}}
\def\llbracket{[\hspace{-.10em} [ }
\def\rrbracket{ ] \hspace{-.10em}]}
\def\cq{$\hfill \square$}
\def \un{\underline}

\def\m{{\cal M}}
\def\tree{{\cal T}_{\ov{\bf e}}}
\def\n{{\cal N}}

\def\u{{\cal U}}

\def\t{{\cal T}}

\def\W{{\cal W}}
\def\w{{\rm w}}
\def\mm{{\bf m}}
\def\eg{{\bf e}}

\def\be{\begin{equation}}
\def\ee{\end{equation}}
\def\ba{\begin{eqnarray*}}
\def\ea{\end{eqnarray*}}
\def\ov{\overline}
\def\wh{\widehat}
\def\wt{\widetilde}

\def\la{\longrightarrow}

\def\proof{\vskip 3mm \noindent{\bf Proof:}\hskip10pt}

\newtheorem{theorem}{Theorem}[section]
\newtheorem{lemma}[theorem]{Lemma}
\newtheorem{proposition}[theorem]{Proposition}

\newtheorem{definition}{Definition}[section]

\begin{document}

\title{ \bf THE TOPOLOGICAL STRUCTURE\\
OF SCALING LIMITS\\
OF LARGE PLANAR MAPS}
\author{
Jean-Fran\c cois {\sc Le Gall}\footnote{DMA-ENS, 45 rue d'Ulm, 75005 Paris, France --- e-mail:
legall@dma.ens.fr\ , fax: (33) 1 44 32 20 80}{}\\ {\small Ecole normale sup\'erieure de Paris}}
\date{\small\today}
\maketitle

\begin{abstract}
We discuss scaling limits of large bipartite planar maps.
If $p\geq 2$ is a fixed integer, we consider, for every integer $n\geq 2$, a random planar
map $M_n$ which is uniformly distributed over the set of all rooted $2p$-angulations
with $n$ faces. Then, at least along a suitable subsequence, the
metric space consisting of the set of vertices of $M_n$, equipped with the graph distance rescaled by the
factor
$n^{-1/4}$, converges in distribution as $n\to\infty$ towards a limiting random compact
metric space, in the sense of the Gromov-Hausdorff distance. We prove that the topology of the limiting
space is uniquely determined independently of $p$ and of the subsequence, and that 
this space can be obtained as the quotient of the
Continuum Random Tree for an equivalence relation which is defined from Brownian labels attached to
the vertices. We also verify that the Hausdorff dimension of the 
limit is almost surely equal to $4$.
\end{abstract}

\section{Introduction}

The main purpose of the present work is to investigate
continuous limits of rescaled planar maps. We concentrate
on bipartite planar maps, which are known to be in
one-to-one correspondence with certain labeled trees
called mobiles (Bouttier, Di Francesco, Guitter \cite{BDG}). 
In view of the
correspondence between maps and mobiles, it seems plausible that scaling limits of large 
bipartite planar maps
can be described in terms of continuous random trees. This 
idea already appeared in the pioneering work of Chassaing
and Schaeffer \cite{CS}, and was then developed by 
Marckert and Mokkadem \cite{MaMo}, who
defined and studied the so-called Brownian map. It was argued in \cite{MaMo} that
the Brownian map is in some weak sense the limit of rescaled uniformly
distributed random quadrangulations of the plane (see also Marckert and Miermont
\cite{MaMi} for recent work along the same lines). The point of view
of the present paper is however different from the one in \cite{MaMo}
or in \cite{MaMi}. For every given planar map $M$, we equip the
set ${\bf m}$ of its vertices with the graph distance, and our aim is to study
the resulting compact metric space when the number of faces 
of the map tends to infinity.  Assuming that the 
map $M$ is chosen uniformly over the set of all rooted $2p$-angulations with 
$n$ faces, we discuss the convergence in distribution when $n$ tends to infinity of the 
associated random metric spaces, rescaled with the factor $n^{-1/4}$,
in the sense of the Gromov-Hausdorff distance between compact metric spaces
(see e.g. Chapter 7 of \cite{BBI}, and subsection 2.3 below, for the definition of the Gromov-Hausdorff
distance). This is in contrast with \cite{MaMo}, which does not consider the limiting 
behavior of distances between two points other than the root vertex.

Before we describe our main results in a more precise way, we need to
set some definitions. Recall that a planar map is a proper embedding, without edge
crossings, of a connected graph in the two-dimensional sphere. Loops and multiple edges
are a priori allowed. The faces of the map are the connected components of the complement
of the union of edges. A planar map is rooted 
if it has a distinguished oriented edge called the root edge,
whose  origin is called the root vertex. The set of vertices will always be equipped
with the graph distance: If $a$ and $a'$ are two vertices, $d_{gr}(a,a')$ is the
minimal number of edges on a path from $a$ to $a'$. Two rooted planar maps
are said to be equivalent if the second one is the image of the first one
under an orientation-preserving
homeomorphism of the sphere, which also preserves the root edges.
From now on we deal only with equivalence classes of rooted planar maps. 
Given an integer $p\geq 2$, a
$2p$-angulation is a planar map where each face has degree $2p$, that is $2p$ adjacent
edges (one should count edge sides, so that if an edge lies entirely inside a face
it is counted twice). We denote by $\m^p_n$ the set of all
rooted $2p$-angulations with $n$ faces.

Let us now discuss the continuous trees that will arise in scaling limits
of planar maps.
We write $\t_\eg$ for the continuum random tree or CRT,
which was introduced and studied by Aldous \cite{Al1}, \cite{Al3}. The CRT 
can be viewed as a random variable taking
values in the space of all rooted compact real trees (see e.g. \cite{survey}, or subsection 2.3
below). It turns out that the CRT 
is the limit in distribution of several (suitably rescaled) classes of discrete trees when the
number of edges tends to infinity. For instance, it is relatively easy to show that if $\tau_n$
is distributed uniformly over the set of all plane trees with $n$ edges, then the vertex set 
of $\tau_n$, viewed as a metric space for the graph distance rescaled by the factor $(2n)^{-1/2}$, will
converge in distribution to the CRT as $n\to\infty$, in the sense of the Gromov-Hausdorff distance. Our
notation
$\t_\eg$ reflects the fact that the CRT can be defined as the real tree coded by a normalized Brownian
excursion $\eg=(\eg_t)_{0\leq t\leq 1}$. This coding, which plays a major role in the
present work, is recalled in subsection 2.3 below. In addition to the usual genealogical order
of the tree, the CRT $\t_\eg$ inherits a lexicographical order from the coding, in a way analogous
to the ordering of (discrete) plane trees from the left to the right. We write $d_\eg$
for the distance on the tree $\t_\eg$ and $\rho$ for the root of $\t_\eg$.

We can assign Brownian labels to the vertices of the CRT. This means that given $\t_\eg$,
we consider a centered Gaussian process $(Z_a)_{a\in\t_\eg}$, such that $Z_\rho=0$  and the variance of
$Z_a-Z_b$ is equal to $d_\eg(a,b)$ for every
$a,b\in\t_\eg$. The pair $(\t_\eg,(Z_a)_{a\in\t_\eg})$ is the probabilistic
object that allows us to describe the continuous limit of random planar maps. We use 
the Brownian labels to define a mapping $D^\circ$ from $\t_\eg\times \t_\eg$
into $\R_+$, via the formula
$$D^\circ(a,b)=Z_a+Z_b - 2\inf_{c\in[a,b]} Z_c$$
where $[a,b]$ denotes the ``lexicographical'' interval between $a$ and $b$. The preceding 
definition is a little informal, since there are two lexicographical intervals between $a$
and $b$, corresponding to the two possible ways of going from $a$ to $b$ around the tree. 
It should be understood that we take the lexicographical interval that minimizes 
the value of $D^\circ(a,b)$ as defined above (see Section 3 below for a more
rigorous presentation). The intuition behind the definition of $D^\circ$ comes from the
discrete picture where each (bipartite) planar map is coded by a labeled tree, in such a way 
that vertices of the map other than the root are in one-to-one correspondence with
vertices of the tree (\cite{BDG}, see subsection 2.1 below). From the properties of this coding, and more
precisely from the way edges of the map are reconstructed from the labels in the tree, one sees that any
two vertices
$a$ and $b$
that satisfy a discrete version of the relation $D^\circ(a,b)=0$ will be connected by
an edge of the map. See subsection 2.1 for more details.

The function $D^\circ$ does not satisfy the triangle inequality, but we may set
$$D^*(a,b)=\inf\left\{\sum_{i=1}^q D^\circ(a_{i-1},a_i)\right\}$$
where the infimum is over all choices of the integer $q\geq 1$
and of the finite sequence $a_0,a_1,\ldots,a_q$
in $\t_{{\bf e}}$ such that $a_0=a$ and $a_q=b$. We then define an equivalence relation
on $\t_{\eg}$ by setting $a\approx b$ if and only if $D^*(a,b)=0$. Although this is not
obvious, it turns out that the latter condition is equivalent to $D^\circ(a,b)=0$, outside a
set of probability zero. 
Moreover one can check that equivalence classes for $\approx$ contain $1$, $2$ or at
most $3$ points, almost surely.
The quotient space $\t_\eg\,/\!\approx$ equipped with the metric $D^*$ is compact.

Let us now come to our main results. For every integer $n\geq 2$, let 
$M_n$ be a random rooted $2p$-angulation uniformly distributed over $\m^p_n$. Denote 
by ${\bf m}_n$ the set of vertices of $M_n$ and by $d_n$ the
graph distance on ${\bf m}_n$. We view $({\bf m}_n,d_n)$ as a random
variable taking values in the space of isometry classes of compact metric spaces.
Recall that the latter space equipped with the Gromov-Hausdorff distance is a Polish space,
as a simple consequence of Gromov's compactness theorem (\cite{BBI}, Theorem 7.4.15).
It can be checked that the sequence of the laws of 
$({\bf m}_n,n^{-1/4}d_n)$ is tight, and so, at least along a subsequence, we
may assume that $({\bf m}_n,n^{-1/4}d_n)$ converges in distribution towards 
a certain random compact metric space. The Skorokhod representation
theorem even allows us to get
an almost sure convergence, at the cost of replacing each map $M_n$
by another random map with same distribution. The principal contribution of the present
work is to identify the limiting compact metric space up to homeomorphism.

Precisely, our main result (Theorem \ref{main}) can be stated as follows.
From any
sequence of integers converging to $+\infty$, we can extract a subsequence and for every $n$
belonging to this subsequence
we can construct a  random $2p$-angulation $M_n$ that is uniformly distributed over $\m^p_n$, in such a way
that we have the almost sure
convergence
\be
\label{mainintro}
\left({\bf m}_n,\Big(\frac{9}{4p(p-1)}\Big)^{1/4}\,n^{-1/4} d_n\right)
\build{\la}_{n\to\infty}^{} (\t_{\eg}\,/\!\approx,D)
\ee
in the sense of the Gromov-Hausdorff distance. Here $D$ is a (random) metric on the quotient space
$\t_{\eg}\,/\!\approx$, such that $D(a,b)\leq D^*(a,b)$ for every $a,b$. The random metric $D$ may 
a priori depend on the choice of the subsequence and on the value of $p$. However,
since $\t_\eg\,/\!\approx$ equipped with the metric $D^*$ is compact and $D\leq D^*$, a standard argument
shows that the metric spaces $(\t_{\eg}\,/\!\approx,D)$ and $(\t_{\eg}\,/\!\approx,D^*)$
are homeomorphic, so that the topological structure of the limit in (\ref{mainintro})
is uniquely determined. In the companion paper \cite{LGP}, we prove that 
$(\t_{\eg}\,/\!\approx,D)$, or equivalently $(\t_{\eg}\,/\!\approx,D^*)$, is a.s.
homeomorphic to the sphere $S^2$. We conjecture that $D=D^*$, and then the
convergence (\ref{mainintro}) would not require the use of a subsequence, and the limit would
not depend on $p$ (the constant $(9/(4p(p-1))^{1/4}$ in (\ref{mainintro})
is relevant mainly because we expect the limit to be independent of $p$). Although we are not able to
prove this, we can derive  enough information about the limiting metric space in (\ref{mainintro}) to prove
that its Hausdorff dimension is equal to $4$ almost surely (Theorem \ref{Hausdim}).

Let us briefly comment on the proof of our main result.
The compactness argument that we use to get the existence of a limit in (\ref{mainintro})
along a suitable subsequence also shows that
this limit can be written as a quotient of the CRT $\t_\eg$ corresponding to a certain random 
pseudo-metric $D$. The point is then to check that $a\approx b$ holds if and
only if $D(a,b)=0$. In other words, the points of the CRT that we
need to identify in order to get the limit in (\ref{mainintro}) 
are given by the equivalence relation
$\approx$, which is defined in terms of $D^*$ or of $D^\circ$. Once we know
that $D\leq D^*$, it is obvious that
$a\approx b$ implies $D(a,b)=0$. The hard core of the proof is thus to check the reverse
implication. The above-mentioned interpretation of the condition $D^\circ(a,b)=0$
in the discrete setting makes
it clear that any two points satisfying this condition must be identified. However,
other pairs of points could conceivably have been identified. 
Roughly speaking, the proof that this is not the
case proceeds as follows.  Given $a$ and $b$ in $\t_\eg$,
we can construct corresponding vertices $a_n$ and $b_n$ in $M_n$ such that
the sequence $(a_n)$ converges to $a$ and the sequence $(b_n)$
converges to $b$, in some suitable sense. The condition $D(a,b)=0$
entails that $d_n(a_n,b_n)=o(n^{1/4})$ as $n\to \infty$. We can then use this estimate
together with some combinatorial considerations and certain delicate properties
of the ``Brownian tree'' $(\t_\eg,(Z_a)_{a\in\t_\eg})$, in order to conclude
that we must have $D^\circ(a,b)=0$.

Let us discuss previous work related to the subject of the present article.
Planar maps were first studied by Tutte \cite{Tu} in connection with his work
on the four colors theorem. Because of their relations with Feynman diagrams,
planar maps soon attracted the attention of specialists of theoretical physics. The 
pioneering papers \cite{tH} and \cite{BIPZ} related enumeration problems for planar
maps with asymptotics of matrix integrals. The interest for random planar
maps in theoretical physics grew significantly when these combinatorial
objects were interpreted as models of random surfaces, especially in
the setting of the theory of quantum gravity (see in particular \cite{Da}
and the book \cite{ADJ}). On the other hand, the idea	 of coding 
planar maps with simpler combinatorial objects such as labeled trees
appeared in Cori and Vauquelin \cite{CV} and was much developed in Schaeffer's thesis
\cite{Sc}. In the present work, we use a version of the bijections between maps and trees 
that was obtained in the recent paper of Bouttier, Di Francesco and Guitter \cite{BDG}.
See Bouttier's thesis \cite{Bo} and the references therein for applications of these bijections to
the statistical physics of random surfaces.
Other applications in the spirit of
the present work can be found in the recent papers \cite{CS}, \cite{MaMi} and \cite{MaMo}
that were mentioned earlier. Note in particular that the random metric space 
$(\t_\eg\,/\!\approx,D^*)$ that is discussed above is essentially equivalent
to the Brownian map of \cite{MaMo}, although the presentation there is 
different. See also \cite{An}, \cite{AS}, \cite{CD} and \cite{Kr} for various results
about random infinite planar triangulations and quadrangulations and their asymptotic
properties.

The paper is organized as follows. Section 2 gives a number of preliminaries
concerning bijections between maps and trees, the coding of real trees and 
the construction of the Brownian tree $(\t_\eg,(Z_a)_{a\in\t_\eg})$. We also
state three important lemmas about the Brownian tree. Section 3 contains our main results. 
The presentation is
slightly different (although equivalent) from the one that is given above, 
because we prefer to argue with the tree $\t_\eg$ re-rooted at the vertex
with the minimal label, and the labels $Z_a$ shifted accordingly so that the label
of the root is still zero. Indeed, it is the genealogical structure of this
re-rooted tree that plays a major role in our approach. Section 4 is devoted to the main step
of our arguments, that is the proof that $D(a,b)=0$ implies $D^\circ(a,b)=0$. 
Section 5 gives the proof of three technical lemmas that were stated in
Section 2. The proofs of these lemmas depend on some rather intricate properties
of Brownian trees, which we found convenient to derive using the path-valued
process called the Brownian snake \cite{Zu}. In order to make most of the paper
accessible to the reader who is unfamiliar with the Brownian snake, we have preferred
to postpone these proofs to Section 5. At last, Section 6 contains the calculation
of the Hausdorff dimension of the limiting metric space.

As a final remark, it is very plausible that our results
can be extended to the more general setting of Boltzmann distributions
on bipartite maps, which is considered in \cite{MaMi}
and in \cite{We}. We have chosen to concentrate on the particular
case of uniform $2p$-angulations for the sake of simplicity and
to keep the present work to a reasonable size.

\medskip
\noindent{\bf Acknowledgments}. I am indebted to Gr\'egory Miermont for a number of very
stimulating discussions. I also thank Fr\'ed\'eric Paulin for several useful
conversations and helpful comments, and Oded Schramm for his remarks on a preliminary version of this work.

\section{Preliminaries}

\subsection{Planar maps and the Bouttier-Di Francesco-Guitter bijection}

Recall that we have fixed an integer $p\geq 2$ and that $\m^p_n$ denotes the set of all
rooted $2p$-angulations with $n$ faces. We start this section with a precise
description of the Bouttier-Di Francesco-Guitter bijection between 
$\m^p_n$ and the set of all $p$-mobiles with $n$ black vertices.

We use the standard formalism for plane trees as found in \cite{Neveu}
for instance. Let
$${\cal U}=\bigcup_{n=0}^\infty \N^n  $$
where $\N=\{1,2,\ldots\}$ and by convention $\N^0=\{\varnothing\}$.
The generation of $u=(u_1,\ldots,u_n)\in\N^n$ is $|u|=n$. If
$u=(u_1,\ldots u_m)$ and 
$v=(v_1,\ldots, v_n)$ belong to $\cal U$, $uv=(u_1,\ldots u_m,v_1,\ldots ,v_n)$
denotes the concatenation of $u$ and $v$. In particular $u\varnothing=\varnothing u=u$.
If $v$ is of the form $v=uj$ for $u\in{\cal U}$ and $j\in\N$, 
we say that $u$ is the {\it father} of $v$, or that $v$ is a {\it child} of $u$.
More generally, if $v$ is of the form $v=uw$ for $u,w\in{\cal U}$, we say that 
$u$ is an {\it ancestor} of $v$, or that $v$ is a {\it descendant} of $u$.

A plane tree $\tau$ is a finite subset of
$\cal U$ such that:
\begin{description}
\item{(i)} $\varnothing\in \tau$.

\item{(ii)} If $v\in \tau$ and $v\ne\varnothing$, the father of
$u$ belongs to $\tau$. 

\item{(iii)} For every $u\in\tau$, there exists an integer $k_u(\tau)\geq 0$
such that $uj\in\tau$ if and only if $1\leq j\leq k_u(\tau)$.
\end{description}
A $p$-tree is a plane tree $\tau$ that satisfies the 
following additional property:
\begin{description}
\item{(iv)} For every $u\in\tau$ such that $|u|$ is odd, $k_u(\tau)=p-1$.
\end{description}

If $\tau$
is a $p$-tree, vertices $u$ of $\tau$ such that $|u|$ is even are called
white vertices, and vertices of $u$ such that $|u|$ is odd are called
black vertices. We denote by $\tau^\circ$ the set of all white vertices
of $\tau$ and by $\tau^\bullet$ the set of all black vertices.
See the left side of Fig.1 for an example of a $3$-tree.

\begin{center}
\begin{picture}(450,140)
\linethickness{1.5pt}
\put(160,0){\circle{12}}
\put(146,-14){$\varnothing$}
\put(155,3){\thinlines\line(-2,1){56}}

\put(165,3){\thinlines\line(2,1){56}}
\put(220,30){\circle*{12}}
\put(228,17){$2$}
\put(245,55){\circle{12}}
\put(223,33){\line(1,1){18}}
\put(195,55){\circle{12}}
\put(180,42){$21$}
\put(217,33){\line(-1,1){18}}
\put(252,42){$22$}

\put(100,30){\circle*{12}}
\put(90,17){$1$}
\put(125,55){\circle{12}}
\put(132,42){$12$}
\put(103,33){\line(1,1){18}}
\put(75,55){\circle{12}}
\put(61,42){$11$}
\put(97,33){\line(-1,1){18}}
\put(75,85){\circle*{12}}
\put(55,72){$111$}
\put(75,61){\thinlines\line(0,1){20}}
\put(100,110){\circle{12}}
\put(78,88){\line(1,1){18}}
\put(50,110){\circle{12}}
\put(23,97){$1111$}
\put(72,88){\line(-1,1){18}}
\put(106,97){$1112$}

\thinlines \put(300,0){\vector(1,0){140}}
\put(300,0){\circle*{3}}
\thinlines \put(300,0){\vector(0,1){130}}
\put(314,56){\circle*{3}}
\thicklines \put(300,0){\line(1,4){14}}
\thicklines \put(314,56){\line(1,4){14}}
\thicklines \put(328,112){\line(1,0){14}}
\put(328,112){\circle*{3}}
\thicklines \put(342,112){\line(1,-4){14}}
\put(342,112){\circle*{3}}
\thicklines \put(356,56){\line(1,0){14}}
\put(356,56){\circle*{3}}
\thicklines \put(370,56){\line(1,-4){14}}
\put(370,56){\circle*{3}}
\thicklines \put(384,0){\line(1,4){14}}
\put(384,0){\circle*{3}}
\thicklines \put(398,56){\line(1,0){14}}
\put(398,56){\circle*{3}}
\thicklines \put(412,56){\line(1,-4){14}}
\put(412,56){\circle*{3}}
\put(426,0){\circle*{3}}

\put(440,-10){$i$}
\put(282,125){$C^{\tau^\circ}_i$}

\put(422,-10){$pn$}

\put(314,0){\line(0,1){4}}
\put(312,-10){$1$}

\put(293,52){$1$}
\put(300,55){\line(1,0){4}}

\end{picture}

\vspace{8mm}

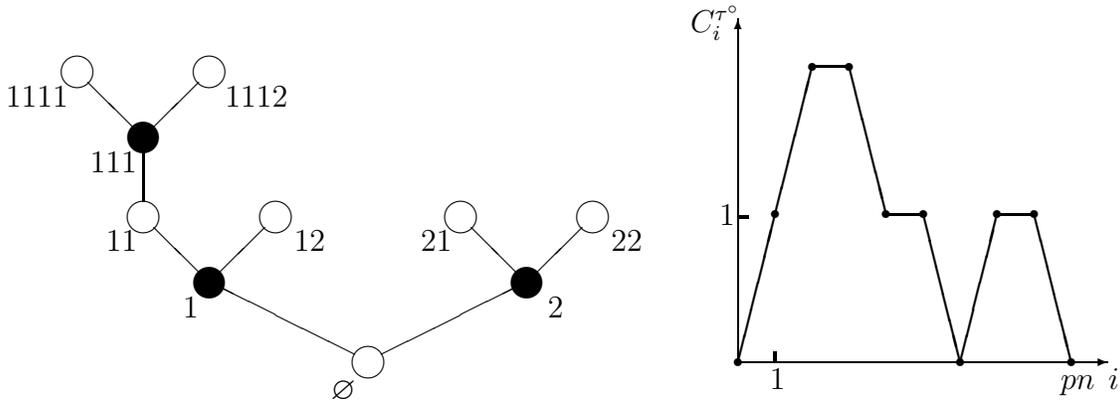
Figure 1. A $3$-tree $\tau$ and the associated 
contour function $C^{\tau^\circ}$ of $\tau^\circ$.

\end{center}

A (rooted) $p$-mobile is a pair $\theta=(\tau,(\ell_u)_{u\in\tau^\circ})$
that consists of a $p$-tree $\tau$ and a collection of integer labels
attached to the white vertices of
$\tau$, such that the following properties hold:
\begin{description}
\item{(a)} $\ell_\varnothing=1$ and $\ell_u\geq 1$ for each $u\in\tau^\circ$.
\item{(b)} Let $u\in \tau^\bullet$, let $u_{(0)}$ be the father of $u$ and let
$u_{(j)}=uj$ for every $1\leq j\leq p-1$. Then for every $j\in\{0,1,\ldots,p-1\}$,
$\ell_{u_{(j+1)}}\geq \ell_{u_{(j)}}-1$, where by convention $u_{(p)}=u_{(0)}$.
\end{description}

\begin{center}
\begin{picture}(400,170)
\linethickness{1.5pt}
\put(140,0){\circle{12}}
\put(137,-4){$1$}
\put(135,3){\line(-4,3){34}}
\put(145,3){\line(4,3){34}}

\put(180,30){\circle*{8}}
\put(205,55){\circle{12}}
\put(183,33){\line(1,1){18}}
\put(155,55){\circle{12}}
\put(152,51){$1$}
\put(177,33){\line(-1,1){18}}
\put(202,51){$2$}

\put(100,30){\circle*{8}}
\put(125,55){\circle{12}}
\put(103,33){\line(1,1){18}}
\put(75,55){\circle{12}}
\put(72,51){$3$}
\put(97,33){\line(-1,1){18}}
\put(122,51){$2$}
\put(75,85){\circle*{8}}
\put(75,61){\thinlines\line(0,1){20}}
\put(100,110){\circle{12}}
\put(78,88){\line(1,1){18}}
\put(50,110){\circle{12}}
\put(47,106){$4$}
\put(72,88){\line(-1,1){18}}
\put(97,106){$3$}
\put(46,114){\line(-1,1){20}}
\put(25,135){\circle*{8}}
\put(104,114){\line(1,1){20}}
\put(125,135){\circle*{8}}
\put(150,160){\circle{12}}
\put(128,138){\line(1,1){18}}
\put(100,160){\circle{12}}
\put(97,156){$2$}
\put(122,138){\line(-1,1){18}}
\put(147,156){$1$}
\put(50,160){\circle{12}}
\put(28,138){\line(1,1){18}}
\put(0,160){\circle{12}}
\put(-3,156){$3$}
\put(22,138){\line(-1,1){18}}
\put(47,156){$2$}

\thinlines \put(250,2){\vector(1,0){180}}
\thinlines \put(250,2){\vector(0,1){160}}

\put(250,35){\circle*{3}}
\put(250,35){\line(1,6){11}}
\put(261,101){\circle*{3}}
\put(261,101){\line(1,3){11}}
\put(272,134){\circle*{3}}
\put(272,134){\line(1,-3){11}}
\put(283,101){\circle*{3}}
\put(283,101){\line(1,-3){11}}
\put(294,68){\circle*{3}}
\put(294,68){\line(1,6){11}}
\put(305,134){\circle*{3}}
\put(305,134){\line(1,-3){11}}
\put(316,101){\circle*{3}}
\put(316,101){\line(1,-3){11}}
\put(327,68){\circle*{3}}
\put(327,68){\line(1,-3){11}}
\put(338,35){\circle*{3}}
\put(338,35){\line(1,6){11}}
\put(349,101){\circle*{3}}
\put(349,101){\line(1,0){11}}
\put(360,101){\circle*{3}}
\put(360,101){\line(1,-3){11}}
\put(371,68){\circle*{3}}
\put(371,68){\line(1,-3){11}}
\put(382,35){\circle*{3}}
\put(382,35){\line(1,3){11}}
\put(393,68){\circle*{3}}
\put(393,68){\line(1,0){11}}
\put(404,68){\circle*{3}}
\put(404,68){\line(1,-3){11}}
\put(415,35){\circle*{3}}

\put(242,32){$1$}
\put(242,65){$2$}
\put(242,98){$3$}
\put(242,131){$4$}
\put(250,68){\line(1,0){3}}
\put(250,101){\line(1,0){3}}
\put(250,134){\line(1,0){3}}

\put(261,2){\line(0,1){3}}
\put(259,-8){$1$}
\put(415,2){\line(0,1){3}}
\put(409,-8){$pn$}

\put(428,-8){$i$}
\put(233,155){$V^\theta_i$}

\end{picture}

\vspace{4mm}

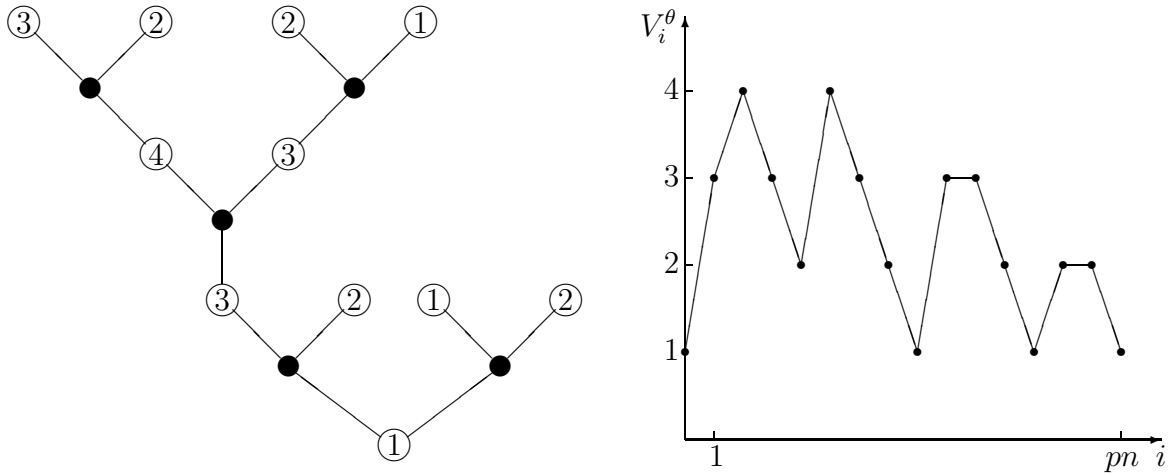
Figure 2. A $3$-mobile $\theta$ with $5$ black vertices
and the associated spatial contour function.

\end{center}

The left side of Fig.2 gives an example of a $p$-mobile with $p=3$. The numbers appearing 
inside the circles representing white vertices are the labels assigned to
these vertices. Condition (b) above means that if one lists the white vertices
adjacent to a given black
vertex in clockwise order, the labels of these vertices can decrease by at most one
at each step.

We will now describe the Bouttier-Di Francesco-Guitter bijection between 
$\m^p_n$ and the set of all $p$-mobiles with $n$ black vertices.
This bijection can be found in Section 2 of \cite{BDG} in the more general
setting of bipartite planar maps. Also \cite{BDG} deals with pointed
planar maps rather than with rooted planar maps. It is however easy to verify that
the results described below are simple consequences of \cite{BDG}.

Let $\tau$ be a $p$-tree with $n$ black vertices and let $k=\#\tau -1=pn$. The
search-depth sequence of $\tau$ is the sequence $u_0,u_1,\ldots,u_{2k}$ of vertices
of $\tau$ which is obtained by induction as follows. First $u_0=\varnothing$, and
then for every $i\in\{0,\ldots,2k-1\}$, $u_{i+1}$ is either the first child of
$u_i$ that has not yet appeared in the sequence $u_0,\ldots,u_i$, or the father
of $u_i$ if all children of $u_i$ already appear in the sequence $u_0,\ldots,u_i$.
It is easy to verify that $u_{2k}=\varnothing$ and that all vertices of $\tau$
appear in the sequence $u_0,u_1,\ldots,u_{2k}$ (of course some of them 
appear more than once).

It is immediate to see that vertices $u_i$ are white when $i$ is even
and black when $i$ is odd.
The search-depth sequence of $\tau^\circ$ is by definition the sequence 
$v_0,\ldots,v_k$ defined by $v_i=u_{2i}$ for every $i\in\{0,1,\ldots,k\}$.

Now let $(\tau,(\ell_u)_{u\in\tau^\circ})$ be a $p$-mobile 
with $n$ black vertices. Denote by $v_0,v_1,\ldots,v_{pn}$ the search-depth
sequence of $\tau^\circ$. Suppose that the tree $\tau_n$ is drawn 
in the plane as pictured on Fig.3 and add an extra vertex $\partial$.
We associate with $(\tau,(\ell_u)_{u\in\tau^\circ})$ a $2p$-angulation $M$
with $n$ faces,
whose set of vertices is
$$\tau^\circ \cup\{\partial\}$$
and whose edges are obtained by the following device: For every
$i\in\{0,1,\ldots,pn-1\}$,
\begin{description}
\item{$\bullet$} if $\ell_{v_i}=1$, draw an edge between $v_i$ and $\partial$\ ;
\item{$\bullet$} if $\ell_{v_i}\geq 2$, draw an edge between $v_i$ and the first
vertex in the sequence $v_{i+1},\ldots,v_{pn}$ whose label is 
$\ell_{v_i}-1$ (this vertex will be called a {\it successor} of $v_i$  -- note that a given 
vertex $v$ can appear several times in the search-depth sequence and so 
may have several different successors).

\end{description}

Notice that $\ell_{v_{pn}}=\ell_\varnothing=1$ and that condition (b) 
in the definition of a $p$-tree
entails that
$\ell_{v_{i+1}}\geq \ell_{v_i}-1$ for every $i\in\{0,1,\ldots,pn-1\}$. 
This ensures that whenever $\ell_{v_i}\geq 2$ there is at least one
vertex among $v_{i+1},v_{i+2},\ldots,v_{pn}$ with label $\ell_{v_i}-1$.
The construction can be made in such a way
that edges do not intersect: See
Section 2 of \cite{BDG}. The resulting planar
graph $M$ is a $2p$-angulation, which is rooted at the oriented edge 
between $\partial$ and $v_0=\varnothing$, corresponding 
to $i=0$ in the previous construction. Each black
vertex of
$\tau$ is associated with a face of the map $M$. Furthermore the graph 
distance in $M$ between the root vertex $\partial$ and another vertex
$u\in\tau^\circ$ is equal to $\ell_u$. See Fig.3 for the $6$-angulation
associated with the $3$-mobile of Fig.2.

It follows from \cite{BDG} that the preceding construction yields
a bijection between the set $\T^p_n$ of all $p$-mobiles with $n$
black vertices and the set $\m^p_n$.

\begin{center}
\begin{picture}(400,220)
\linethickness{1.5pt}
\put(260,0){\circle{12}}
\put(257,-4){$1$}
\put(255,3){\line(-2,1){52}}
\put(265,3){\line(4,3){34}}

\put(300,30){\circle*{8}}
\put(325,55){\circle{12}}
\put(303,33){\line(1,1){18}}
\put(275,55){\circle{12}}
\put(272,51){$1$}
\put(297,33){\line(-1,1){18}}
\put(322,51){$2$}

\put(200,30){\circle*{8}}
\put(225,55){\circle{12}}
\put(203,33){\line(1,1){18}}
\put(175,55){\circle{12}}
\put(172,51){$3$}
\put(197,33){\line(-1,1){18}}
\put(222,51){$2$}
\put(175,85){\circle*{8}}
\put(175,61){\thinlines\line(0,1){20}}
\put(200,110){\circle{12}}
\put(178,88){\line(1,1){18}}
\put(150,110){\circle{12}}
\put(147,106){$4$}
\put(172,88){\line(-1,1){18}}
\put(197,106){$3$}
\put(146,114){\line(-1,1){20}}
\put(125,135){\circle*{8}}
\put(204,114){\line(1,1){20}}
\put(225,135){\circle*{8}}
\put(250,160){\circle{12}}
\put(228,138){\line(1,1){18}}
\put(200,160){\circle{12}}
\put(197,156){$2$}
\put(222,138){\line(-1,1){18}}
\put(247,156){$1$}
\put(150,160){\circle{12}}
\put(128,138){\line(1,1){18}}
\put(100,160){\circle{12}}
\put(97,156){$3$}
\put(122,138){\line(-1,1){18}}
\put(147,156){$2$}
\bezier{800}(168,55)(-60,260)(146,165)
\bezier{50}(106,160)(125,170)(143,160)
\bezier{100}(143,110)(100,132)(100,154)
\bezier{50}(156,110)(175,120)(193,110)
\bezier{50}(200,116)(190,135)(200,154)
\bezier{100}(153,165)(200,210)(250,166)
\bezier{50}(206,160)(225,170)(243,160)
\bezier{100}(206,110)(230,85)(225,61)
\bezier{50}(181,55)(200,65)(218,55)
\bezier{80}(225,49)(230,22)(256,5)

\put(300,140){\circle{12}}
\put(297,136){$\partial$}
\bezier{200}(260,6)(250,110)(294,140)
\bezier{200}(256,160)(280,182)(300,146)
\bezier{200}(275,61)(300,100)(300,134)
\bezier{150}(325,49)(320,30)(266,0)
\bezier{200}(265,-3)(400,20)(306,140)

\end{picture}

\vspace{6mm}
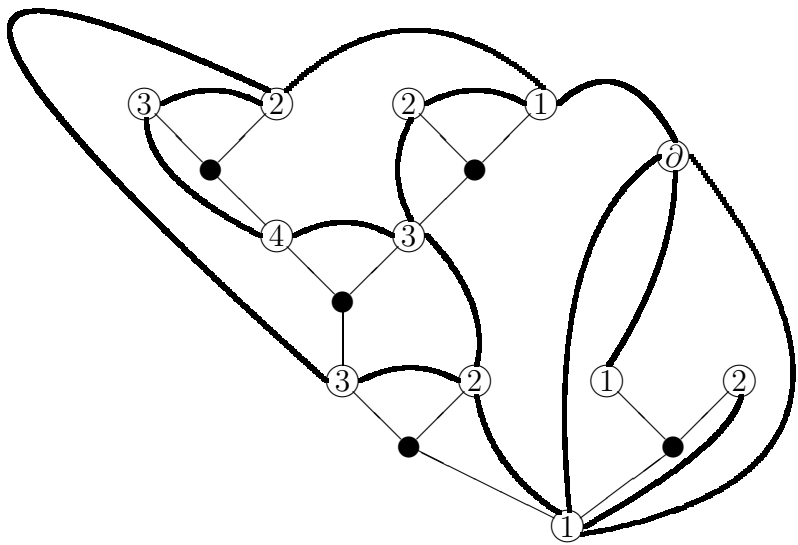
Figure 3. The Bouttier-Di Francesco-Guitter bijection: A rooted $3$-mobile with $5$ black vertices
and the associated
$6$-angulation with $5$ faces

\end{center}

\subsection{Genealogical structure of maps}

Let $\theta=(\tau,(\ell_u)_{u\in\tau^\circ})$ be a $p$-mobile
with $n$ black vertices. The set 
$\tau^\circ$ of white vertices can also be viewed as a graph, by declaring that 
there is an edge between $u$ and $v$ if and only if $u$ is the grandfather
of $v$ (that is, there exist $j$ and $k\in\N$ such that $v=ujk$) or conversely
$v$ is the grandfather of $u$. Obviously $\tau^\circ$ is a tree in
the graph-theoretic sense. If $u,v\in\tau^\circ$, we then denote by
$\llbracket u,v\rrbracket$ the set of points of $\tau^\circ$ that lie on
the unique shortest path from $u$ to $v$ in $\tau^\circ$. As
usual,
$\rrbracket u,v\llbracket=
\llbracket u,v\rrbracket\backslash\{u,v\}$. We also denote by $u\wedge v$
the ``most recent common ancestor'' of $u$ and $v$ in $\tau^\circ$, which may be
defined by $\llbracket \varnothing,u\wedge v\rrbracket=
\llbracket \varnothing,u\rrbracket\cap \llbracket \varnothing, v\rrbracket$.
Notice that $u\wedge v$ is not necessarily the most recent ancestor of
$u$ and $v$ in the tree $\tau$.

We denote
by $\prec$ the genealogical relation on $\tau^\circ$: $u\prec v$ if and only
if $u$ is an ancestor of $v$ (in the tree $\tau$). We use $u\leq v$ for the 
lexicographical order on $\tau^\circ$. As usual $u<v$ if and
only if $u\leq v$ and $u\ne v$. It will also be
convenient to introduce a ``reverse'' lexicographical order denoted by $\ll$. This
is the total order on $\tau^\circ$ defined as follows. If neither of the
relations
$u\prec v$ and $v\prec u$ holds, then $u\ll v$ if and only if $u\leq v$. On
the other hand, if $u\prec v$, then $v\ll u$ (although $u \leq v$).

Let $v_0,v_1,\ldots,v_{pn}$ be the search-depth sequence of $\tau^\circ$, as defined in
the preceding subsection. If $x,y\in\tau^\circ$, the condition $x\leq y$
implies that the first occurence of $x$ in the sequence $v_0,\ldots,v_{pn}$ 
occurs before the first occurence of $y$, and conversely the condition $x\ll y$
implies that the last occurence of $x$ occurs before the last occurence of $y$.
The {\it contour function} of $\tau^\circ$ is the discrete sequence
$C^{\tau^\circ}_0,C^{\tau^\circ}_1,\ldots,C^{\tau^\circ}_{pn}$ defined by 
$$C^{\tau^\circ}_i=\frac{1}{2}\,|v_i|\ ,\hbox{ for every }0\leq i\leq pn.$$
See Fig.1 for an example with $p=n=3$.
It is easy to verify that the contour function determines $\tau^\circ$, which
in turn determines the $p$-tree $\tau$ uniquely. We will also
use the {\it spatial contour function} of $\theta=(\tau,(\ell_u)_{u\in\tau^\circ})$, which
is the discrete sequence $(V^\theta_0,V^\theta_1,\ldots,V^\theta_{pn})$ defined by
$$V^\theta_i=\ell_{v_i}\ ,\hbox{ for every }0\leq i\leq pn.$$
From property (b) of the labels and the definition of the search-depth sequence, it
is clear that $V^\theta_{i+1}\geq V^\theta_i-1$ for every $0\leq i\leq pn-1$
(cf Fig.2). This fact will
be used many times below.

The pair $(C^{\tau^\circ},V^\theta)$ determines $\theta$ uniquely. For our 
purposes it will sometimes be convenient to view 
$C^{\tau^\circ}$ or $V^\theta$ as functions of the continuous parameter
$t\in[0,pn]$, simply by interpolating linearly on the intervals
$[i-1,i]$, $1\leq i\leq n$ (as it is suggested by Figs 1 and 2).

Let $[pn]$ stand for the set $\{0,1,\ldots,pn\}$. Define an equivalence
relation $\sim$ on $[pn]$ by setting $i\sim j$
if and only if $v_i=v_j$. The quotient space $[pn]/\sim$
is then obviously identified with $\tau^\circ$. This
identification plays an important role throughout this work. 
If $i\leq j$, the
relation $i\sim j$ implies
$$\inf_{i\leq k\leq j} C^{\tau^\circ}_k=C^{\tau^\circ}_i=C^{\tau^\circ}_j.$$
The converse is not true (except if $p=2$) but 
the conditions $j>i+1$, $C^{\tau^\circ}_i=C^{\tau^\circ}_j$ and
$$C^{\tau^\circ}_k>C^{\tau^\circ}_i\ ,\hbox{ for every }
k\in]i,j[\cap \Z$$
imply that $i\sim j$. Similarly, if $i<j$, the condition $v_i\prec v_j$ implies
$$\inf_{i\leq k\leq j} C^{\tau^\circ}_k=C^{\tau^\circ}_i.$$
The converse is not true, but the condition
$$\inf_{i< k\leq j} C^{\tau^\circ}_k>C^{\tau^\circ}_i$$
forces $v_i\prec v_j$.

Let $u,v\in\tau^\circ$ with $u\prec v$, and let $w\in\rrbracket u,v\llbracket$.
The set 
$$\tau^\circ_{(v,w)}:=\{x\in\tau^\circ: x\wedge v=w\hbox{ and }x\leq v\}$$
is called the subtree from the left side of $\llbracket u,v\rrbracket$
with root $w$. Similarly, the set 
$$\widetilde\tau^\circ_{(v,w)}
:=\{x\in\tau^\circ: x\wedge v=w\hbox{ and }v\ll x\}$$
is called the subtree from the right side of $\llbracket u,v\rrbracket$
with root $w$. Let $j\in[pn]$ be such that $v_j=v$, and
set
\ba
&&k=\inf\{i\in\{0,1,\ldots,j\}:v_i=w\},\\
&&k'=\sup\{i\in\{0,1,\ldots,j\}:v_i=w\}.
\ea
Then $\tau^\circ_{(v,w)}$ exactly consists of the vertices $v_i$
for $k\leq i\leq k'$: We will say that
$[k,k']\cap\Z$ is the interval coding $\tau^\circ_{(v,w)}$. Similar remarks apply to
$\widetilde\tau^\circ_{(v,w)}$.

Recall from the preceding subsection that the $p$-mobile
$(\tau,(\ell_u)_{u\in\tau^\circ})$ corresponds to 
a $2p$-angulation $M$ via the Bouttier-Di Francesco-Guitter bijection.
Through this correspondence, vertices of $M$ (with the exception of the root vertex $\partial$)
are identified with elements of $\tau^\circ$. From now on, we systematically do
this identification.
Let $d_M$ stand for the graph distance on the set of vertices of $M$. A geodesic path in $M$ is
a discrete path $\gamma=(\gamma(i),0\leq i\leq k)$ in $M$ such that 
$d_M(\gamma(i),\gamma(j))=|i-j|$ for every $i,j\in\{0,\ldots,k\}$.

The following lemma plays an important role in our proofs.

\begin{lemma}
\label{combi}
Let $\gamma=(\gamma(i),0\leq i\leq k)$ be a geodesic path in $M$ which does not
visit the root vertex $\partial$. Let 
$u=\gamma(0)$ be the starting point of the path $\gamma$ and let $y=\gamma(k)$
be its final point.
Let $\tau_1$ be a subtree from the left
side of $\llbracket \varnothing,y\rrbracket$ (respectively 
from the right
side of $\llbracket \varnothing,y\rrbracket$) with root $w\in\rrbracket
\varnothing,y\llbracket$. Let $v=\gamma(1)$ be the point following $u$ on
the path $\gamma$, and assume that:
\begin{description}
\item{\rm (i)} $v\in\tau_1$ and $v\ne w$.
\item{\rm (ii)} $u\leq w$ (resp. $w\ll u$).
\item{\rm (iii)} For every $i\in\{1,\ldots,k\}$, one has $w\leq \gamma(i)$
(resp.
$\gamma(i)\ll w$).
\end{description}
Then, for any point $x$ of $\tau_1\backslash\{w\}$ such that
\be
\label{combitech1}
\ell_z>\sup_{0\leq i\leq k} \ell_{\gamma(i)},\hbox{ for every }z\in\llbracket w,x
\rrbracket
\ee
one has
$$d_M(x,y)\leq d_M(u,y)+\ell_x-\inf_{0\leq i\leq k} \ell_{\gamma(i)}.$$
\end{lemma}

\proof We only treat the case when $\tau_1$ is a subtree from the left
side of $\llbracket \varnothing,y\rrbracket$. We fix a point 
$x\in\tau_1\backslash\{w\}$ such that (\ref{combitech1}) holds. Denote by $b$ the
first point on the geodesic $\gamma$ such that $x\leq b$. This makes sense
because $x\leq y$ by the definition of subtrees. Also $b\ne u$ because $u\leq w$
and $w<x$. So we can also introduce the point $a$ preceding $b$
on the geodesic $\gamma$.

Let us first assume that $b\ne v$, or equivalently $a\ne u$. Then 
$w\leq a$ by (iii), and $a\leq x$, which forces $a\in \tau_1$. 
On the other hand, assumption (\ref{combitech1}) guarantees that
$a\notin\llbracket w,x
\rrbracket$. Since  $a\leq x\leq b$, it follows that $a$ cannot be an
ancestor of $b$. Any occurence of $a$ in the search-depth sequence of
$\tau^\circ$ thus happens before the first occurence of $b$ in this sequence.
Now notice that $a$ and $b$ are connected by an edge of the
map $M$, and recall the construction of these edges at the end
of the preceding subsection. It follows that $\ell_z\geq \ell_a$ for every vertex
$z$ such that $a\ll z <b$, whereas $\ell_b=\ell_a-1$.
Note that $a\ll x<b$ (the case $a\prec x$ is excluded since $a\in\tau_1$
and $a\notin\llbracket w,x
\rrbracket$, and $x=b$ is impossible by (\ref{combitech1})), so that the previous
sentence applies to $z=x$. Set $q=\ell_x-\ell_b\geq 1$. We let $i_0$ be the first
index
such that
$v_{i_0}=x$, and observe that $x\leq v_i$ for every $i\geq i_0$. We then define $i_1,\ldots,i_q$ by
setting 
$$i_j=\inf\{i\geq i_0:\ell_{v_i}=\ell_x-j\}\hbox{ for every }1\leq j\leq q.$$
By the preceding considerations, we have $\ell_z\geq \ell_a=\ell_x-q+1$
for every $z$ such that $x\leq z<b$. It follows that $v_{i_q}=b$. On the other hand, $v_{i_0}=x$
and $d_M(v_{j},v_{j+1})=1$ for every $0\leq j\leq q-1$, by the construction of
edges in
$M$. We thus get
$$d_M(x,b)\leq q=\ell_x-\ell_b.$$
Finally,
$$d_M(x,y)\leq d_M(b,y)+d_M(x,b)\leq d_M(u,y)+\ell_x-\ell_b,$$
which gives the desired bound.

In the case when $a=u$ and $b=v$, the argument is almost the same. Note that
$u<w$ ($u=w$ is excluded by (\ref{combitech1})) and $x<b=v$ as previously.
The existence of an edge between $u$ and $v$ warrants that $\ell_v=\ell_u-1$
and that $\ell_z\geq \ell_u$ for every vertex $z$ of $\tau_1$
such that $z<v$. In the same way as before, we get $d_M(x,v)\leq \ell_x-\ell_v$
which leads to the desired bound. \cq  

\subsection{Real trees}

We will now discuss the continuous trees that are scaling limits of
our discrete plane trees. We start with a basic definition.

\begin{definition}
A metric space $(\t,d)$ is a real tree if the following two
properties hold for every $a,b\in \t$.

\begin{description}
\item{\rm(i)} There is a unique
isometric map
$f_{a,b}$ from $[0,d(a,b)]$ into $\t$ such
that $f_{a,b}(0)=a$ and $f_{a,b}(
d(a,b))=b$.
\item{\rm(ii)} If $q$ is a continuous injective map from $[0,1]$ into
$\t$, such that $q(0)=a$ and $q(1)=b$, we have
$$q([0,1])=f_{a,b}([0,d(a,b)]).$$
\end{description}

\noindent A rooted real tree is a real tree $(\t,d)$
with a distinguished vertex $\rho=\rho(\t)$ called the root.
\end{definition}

In what follows, real trees will always be rooted and compact, even if this
is not mentioned explicitly. 

Let us consider a rooted real tree $(\t,d)$.
The range of the mapping $f_{a,b}$ in (i) is denoted by
$\llbracket a,b\rrbracket$ (this is the line segment between $a$
and $b$ in the tree), and we also use the obvious notation 
$\rrbracket a,b\llbracket$. 
In particular, for every $a\in \t$, $\llbracket \rho,a\rrbracket$ is the path 
going from the root to $a$, which we will interpret as the ancestral
line of vertex $a$. More precisely we can define a partial order on the
tree, called the genealogical order, by setting $a\prec b$
if and only if $a\in\llbracket \rho,b
\rrbracket$. 
If $a,b\in\t$, there is a unique $c\in\t$ such that
$\llbracket \rho,a
\rrbracket\cap \llbracket \rho,b
\rrbracket=\llbracket \rho,c
\rrbracket$. We write $c=a\wedge b$ and call $c$ the most recent
common ancestor to $a$ and $b$. The multiplicity of a vertex $a\in\t$ is the number of
connected components of $\t\backslash\{a\}$. In particular, $a$ is
called a leaf if it has multiplicity one.

In a way similar to the discrete case, real trees can be coded by
``contour functions''. If $E$ and $F$ are two topological spaces, we write $C(E,F)$ for the space
of all continuous functions from $E$ into $F$. 
Let $\sigma>0$ and let 
$g\in C([0,\sigma],[0,\infty[)$ be such that $g(0)=g(\sigma)=0$. To
avoid trivialities, we will also assume that $g$ is not identically zero.
For every $s,t\in[0,\sigma]$, we set
$$m_g(s,t)=\inf_{r\in[s\wedge t,s\vee t]}g(r),$$
and
$$d_g(s,t)=g(s)+g(t)-2m_g(s,t).$$
It is easy to verify that $d_g$ is a pseudo-metric 
on $[0,\sigma]$. As usual, we introduce the equivalence
relation
$s\simeq_g t$ if and only if $d_g(s,t)=0$ (or equivalently if and only if $g(s)=g(t)=m_g(s,t)$).
The function $d_g$ induces a distance on the quotient space $\t_g
:=[0,\sigma]\,/\!\simeq_g$, and we keep the
notation $d_g$ for this distance. We denote by
$p_g:[0,\sigma]\longrightarrow
\t_g$ the canonical projection. Clearly $p_g$ is continuous (when
$[0,\sigma]$ is equipped with the Euclidean metric and $\t_g$ with the
metric $d_g$), and therefore $\t_g=p_g([0,\sigma])$ is a compact metric space.

By Theorem 2.1 of \cite{DuLG}, 
the metric space $(\t_g,d_g)$ is a real tree. 
We will always view $(\t_g,d_g)$ as a rooted real tree with root $\rho=p_g(0)=p_g(\sigma)$.
Then, if $s,t\in[0,\sigma]$, the
property $p_g(s)\prec p_g(t)$ holds if and only if $g(s)=m_g(s,t)$.

Let us recall the definition of the Gromov-Hausdorff distance. Let
$(E_1,d_1)$ and $(E_2,d_2)$ be two compact metric spaces. 
The Gromov-Hausdorff distance between $(E_1,d_1)$
and $(E_2,d_2)$ is
$$d_{GH}(E_1,E_2)=\inf\Big(d_{Haus}(\varphi_1(E_1),\varphi_2(E_2))\Big),$$ 
where the infimum is over all isometric embeddings $\varphi_1:E_1\la E$ and
$\varphi_2:E_2\la E$ of $E_1$ and $E_2$ into the same metric
space $(E,d)$, and $d_{Haus}$ stands for the usual Hausdorff distance 
between compact subsets of $E$. 
Then Lemma 2.3 of \cite{DuLG} shows that $\t_g$ depends continuously on $g$,
in the sense that
$$d_{GH}(\t_g,\t_g')\leq 2 \|g-g'\|$$
where $\|g-g'\|$ is the supremum norm of $g-g'$.

In addition to the genealogical order $\prec$, the tree 
$\t_g$ inherits a lexicographical order from the coding 
through the function $g$. Precisely if $a,b\in\t_g$ we
write $a\leq b$ if and only if $s\leq t$, where 
$s$, respectively $t$, is the smallest representative 
of $a$, resp. of $b$, in $[0,\sigma]$.  
We can also introduce a ``reverse'' lexicographical order $\ll$, by replacing smallest by greatest
in the previous sentence. If neither
of the relations $a\prec b$ or $b\prec a$ holds, we have $a\ll b$
if and only if $a\leq b$. On the other hand, if $a\prec b$, we have $b\ll a$.

Let $a,b\in \t_g$. If $a\leq b$, or if $a\ll b$, we define the lexicographical interval 
$[a,b]$ as the image
under the projection
$p_g$ of the minimal interval $[s,t]$ such that $s\leq t$, $p_g(s)=a$ and $p_g(t)=b$. If 
neither of the relations $a\leq b$ or $a\ll b$ holds, then there is no such interval
and we take $[a,b]=\varnothing$.
If $[a,b]$ is nonempty,
then $\llbracket a, b\rrbracket\subset
[a,b]$. Furthermore if $a\prec b$, then both $[a,b]$ and
$[b,a]$ are nonempty, and $[a,b]\cap [b,a]=\llbracket a,b\rrbracket$.

Let $a,b\in \t_g$ with $a\prec b$, and let $c\in\rrbracket a,b\llbracket$.
Suppose that the set
$$\t^1=\{u\in\t_g: u\wedge b =c\hbox{ and }u\leq b\}$$
is not the singleton $\{c\}$. 
Then the set $\t^1$ is called a
subtree from the left side of $\llbracket a,b\rrbracket$ with root $c$
(it is straightforward to verify that $\t^1$ is itself a real tree).
Moreover, if $s=\inf p_g^{-1}(a)$ and $t=\inf p_g^{-1}(b)$, there is
a unique subinterval $[\alpha,\beta]$ of $]s,t[$ such that 
$\t^1=p_g([\alpha,\beta])$, $p_g(\alpha)=p_g(\beta)=c$ and 
$$\alpha=\sup\{r\in[s,t]: g(r)<g(\alpha)\}\ ,\ 
\beta=\sup\{r\in[s,t]: g(r)\leq g(\alpha)\}.$$
We say that $[\alpha,\beta]$ is the coding interval of $\t^1$. In a
similar way we can define subtrees from the right side of $\llbracket
a,b\rrbracket$: $\t^2$ is such a subtree if there exists
$c'\in\rrbracket a,b\llbracket$ such that
$$\t^2=\{u\in\t_g: u\wedge b =c'\hbox{ and }b\ll u\}$$
and $\t^2\ne\{c'\}$.

\subsection{Brownian trees and conditioned Brownian trees}

We first explain how we can assign Brownian labels to 
the vertices of the real tree $(\t_g,d_g)$ defined in
the previous subsection. To this end, we consider the centered real-valued
Gaussian process $(\Gamma_t)_{t\in[0,\sigma]}$ with covariance function
\begin{equation}
\label{covariance}
{\rm cov}(\Gamma_s,\Gamma_t)=m_g(s,t)
\end{equation}
for every $s,t\in[0,\sigma]$ (it is a simple exercise
to check that $m_g(s,t)$ is a covariance function). Note that $\Gamma_0=\Gamma_\sigma=0$ and that the form
of the covariance gives 
$E[(\Gamma_s-\Gamma_t)^2]=d_g(s,t)$. Suppose that $g$ is H\"older
continuous with some exponent $\delta>0$, which will always 
hold in what follows. Then an application of the
classical Kolmogorov lemma shows that the process 
$(\Gamma_t)_{t\in[0,\sigma]}$ has a continuous modification, and from
now on we consider only this modification. We write ${\bf Q}_g$
for the distribution of $(\Gamma_t)_{t\in[0,\sigma]}$, which is
a probability measure on the space $C([0,\sigma],\R)$. 

From the formula $E[(\Gamma_s-\Gamma_t)^2]=d_g(s,t)$ and
a continuity argument, we immediately get that a.s. for every
$s,t\in[0,\sigma]$ such that $s\simeq_g t$, we have $\Gamma_s=\Gamma_t$. 
Therefore we may also view $\Gamma$ as a Gaussian process
indexed by the tree $\t_g$. Indeed, it is natural to interpret
$(\Gamma_a,a\in \t_g)$ as Brownian motion indexed by
$\t_g$ and started from $0$ at the root of $\t_g$. 
Note that formula (\ref{covariance}) may be rewritten in the form
$${\rm cov}(\Gamma_a,\Gamma_b)=d_g(\rho,a\wedge b)$$
for every $a,b\in\t_g$.

We now randomize the coding function $g$.
Let $\eg=(\eg_t)_{t\in[0,1]}$ be the normalized Brownian excursion, and take 
$g=\eg$ and $\sigma=1$
in the previous discussion.
The random real tree $(\t_\eg,d_\eg)$ coded by $\eg$ is the so-called CRT, or
Continuum Random Tree. Using the fact that local minima of Brownian motion
are distinct, one easily checks that points of $\t_\eg$ can have multiplicity at most $3$.

 We then consider
the real-valued  process $(Z_t)_{t\in[0,1]}$ such that conditionally given $\eg$, 
$(Z_t)_{t\in[0,1]}$ has distribution ${\bf Q}_\eg$. As explained
above, we can also view $(Z_t)_{t\in[0,1]}$ as parametrized
by the tree $\t_\eg$, and then interpret $(Z_a)_{a\in\t_\eg}$
as Brownian motion indexed by $\t_\eg$. This interpretation creates some
technical difficulties since $\t_\eg$ is now a random index set -- to
circumvent these difficulties it is often more convenient to view
$Z$ as indexed by $[0,1]$, keeping in mind that $Z_t$ only
depends on the equivalence class of $t$ in $\t_\eg$.
 
In view of our applications it is important 
to consider the pair $(\eg,Z)$
conditioned on the event 
$$Z_t\geq 0\ \hbox{ for every }t\in[0,1].$$
Here some justification is needed for the conditioning,
since the latter event has probability zero. The paper \cite{LGW}
describes several limit procedures that allow one to make sense 
of the previous conditioning. These procedures all lead to
the same limiting pair $(\ov\eg,\ov Z)$ which can be described as 
follows from the original pair $(\eg,Z)$. Set
$$\un Z=\inf_{t\in[0,1]} Z_t$$
and let $s_*$ be the (almost surely) unique time in $[0,1]$ such 
that $Z_{s_*}=\un Z$. The fact that $\un Z$ is attained at a unique
time (\cite{LGW} Proposition 2.5) entails that the vertex $p_\eg(s_*)$ is a 
leaf of the tree $\t_\eg$. For
every $s,t\in[0,1]$, set $s\oplus t=s+t$ if $s+t\leq 1$ and $s\oplus t=s+t-1$ if $s+t>1$. Then, for every
$t\in[0,1]$,
\begin{description}
\item{$\bullet$} 
$\displaystyle{\ov\eg_t=\eg_{s_*}+\eg_{s_*\oplus t}-2\,m_{\eg}(s_*,
{s_*\oplus t})}$;
\item{$\bullet$} $\ov Z_t=Z_{s_*\oplus t} -Z_{s_*}$.
\end{description}
The formula for $\ov Z$ makes it obvious that $\ov Z_t\geq 0$
for every $t\geq 0$, in agreement with the above-mentioned conditioning. The function 
$\ov\eg$ is continuous on $[0,1]$ and such that $\ov\eg(0)=\ov\eg(1)=0$.
Hence the tree $\t_{\ov\eg}$ is well defined, and this tree is
isometrically identified with the tree $\t_\eg$ re-rooted at the
(minimizing) vertex $p_\eg(s_*)$: See Lemma 2.2 in \cite{DuLG}.
Moreover we have $s\simeq_{\ov \eg} t$ if and only if $s_*\oplus s\simeq_{\eg}s_*\oplus t$
and so $\ov Z_t$ only depends on the equivalence class of $t$
in the tree $\t_{\ov\eg}$. Therefore we may and will sometimes
view $\ov Z$ as indexed by vertices of the tree $\t_{\ov\eg}$.

By a well-known property of the Brownian excursion, the law of pair $(\eg_t,Z_t)_{t\in[0,1]}$
is invariant under time reversal, meaning that $(\eg_t,Z_t)_{t\in[0,1]}$
has the same distribution as $(\eg_{1-t},Z_{1-t})_{t\in[0,1]}$. A similar
time-reversal invariance property then holds for the pair
$(\ov\eg_t,\ov Z_t)_{t\in[0,1]}$.
In what follows we use the notation $\rho$ for the root of $\t_\eg$
and $\ov \rho$ for the root of $\t_{\ov\eg}$.

We now state three important lemmas which are key ingredients
of the proofs of our main results. 

\begin{lemma}
\label{increasepoint}
We say that $s\in[0,1[$ is an increase point of the pair
$(\eg,Z)$, respectively of the pair $(\ov\eg,\ov Z)$, if there 
exists $\varepsilon>0$ such that $\eg_t\geq \eg_s$ and
$Z_t\geq Z_s$, resp. $\ov\eg_t\geq \ov\eg_s$ and
$\ov Z_t\geq \ov Z_s$, for every $t\in[s,(s+\varepsilon)\wedge 1]$.
Then a.s. there is no increase point of 
$(\eg,Z)$, and $s=0$ is the only increase point of $(\ov\eg,\ov Z)$.
\end{lemma}

Before stating
the next lemma we need to introduce some additional notation. The
uniform measure $\lambda$ on $\t_\eg$, resp. on $\t_{\ov\eg}$, is the 
image of Lebesgue measure on $[0,1]$ under the canonical projection
$p_{\eg}$, resp. $p_{\ov\eg}$. There is no ambiguity in using the same
notation $\lambda$ for both cases, since it really corresponds to the same
measure when $\t_{\ov\eg}$ is identified to $\t_\eg$ up to re-rooting.
We also let $\cal I$ and $\ov{\cal I}$ be the random measures on $\R$
defined by
$$\langle {\cal I},f\rangle =\int_{\t_\eg} \lambda(da)\,f(Z_a)=\int_0^1 dt\,f(Z_t)
\ ,\ \langle \ov{\cal I},f\rangle =\int_{\t_{\ov\eg}} \lambda(da)\,f(\ov Z_a)=\int_0^1
dt\,f(\ov Z_t).$$
The random measure $\cal I$ is sometimes called (one-dimensional)
ISE. Notice that $\ov{\cal I}$ is supported on $[0,\infty[$ and is just the image of
$\cal I$ under the shift $x\la x-\un Z$. 

\begin{lemma}
\label{estimateISE}
For every $\alpha>0$,
$$\lim_{\varepsilon\to 0} \varepsilon^{-2} P(\ov{\cal I}([0,\varepsilon])\geq
\alpha\varepsilon^2)=0.$$
\end{lemma}

Our last lemma is concerned with values of 
$\ov Z$ over subtrees of $\t_{\ov\eg}$. Roughly speaking it asserts
that, for a given $\beta>0$ and a subtree $\t^1$ with root $c$, if both
$\ov Z_c>\beta$ and the minimum of the values
of $\ov Z$ over $\t^1$ is strictly less than $\beta$, then 
the mass (for the uniform measure $\lambda$) of those vertices $x$ of 
$\t_{\ov\eg}$ with label $\ov Z_x\in[\beta,\beta+\varepsilon]$,
and such that the label of any ancestor of $x$ in $\t^1$ 
is greater than $\beta$, will be of
order at least
$\varepsilon^2$. The precise statement is as follows.

\begin{lemma}
\label{occuptree}
Almost surely, for every $\mu>0$, for every $a\in \t_{\ov\eg}$
and every subtree $\t^1$ from $\llbracket \ov \rho,a\rrbracket$ 
with root $c\in \rrbracket \ov \rho,a\llbracket$,
the condition
$$\inf_{b\in \t^1} \ov Z_b< \ov Z_{c}-\mu$$
implies that
$$\liminf_{\varepsilon\to 0} \varepsilon^{-2}\,\lambda\Big(\Big\{ x\in\t^1: \ov Z_x\leq
\ov Z_c-\mu+\varepsilon
\hbox{ and }\ov Z_y\geq \ov Z_c-\mu+\frac{\varepsilon}{8}
\hbox{ for every }y\in\llbracket c,x\rrbracket\Big\}\Big)
>0.$$
\end{lemma}

Although Lemma \ref{occuptree} is stated in terms of the pair
$(\ov\eg,\ov Z)$, in view of our applications, the proof will show
that this lemma reduces to a similar statement for the pair $(\eg,Z)$.

The proof of the preceding three lemmas depends on some properties of
the path-valued process called the Brownian snake, and recalling these
properties at the present stage would take us too far from our
main concern. For this
reason, we prefer to postpone the proofs to Section 5.

\subsection{Invariance principles}

In this subsection, we recall the basic invariance principles that
relate the discrete labeled trees of subsection 2.1 to the Brownian trees
of subsection 2.4. Recall that the integer $p\geq 2$ is fixed.

Let $\theta_n=(\tau_n,(\ell^n_u)_{u\in\tau^\circ_n})$ be uniformly distributed
over the set $\T^p_n$ of all $p$-mobiles with $n$ black vertices. 
We denote by $C^n=(C^n_t)_{0\leq t\leq pn}$ the contour function of
$\tau_n^\circ$ and by
$V^n=(V^n_t)_{0\leq t\leq pn}$ the spatial contour function of $\theta_n$
(it
is convenient to view $C^n$ and $V^n$ as continuous functions of
$t\in[0,pn]$, as explained in subsection 2.2).
Recall that the  pair $(C^n,V^n)$ determines $\theta_n$. 

\begin{theorem}
\label{invar1}
We have 
\be
\label{basicinvar}
\left({1\over 2}\sqrt{\frac{p}{p-1}}\,n^{-1/2}\,C^n_{pnt},
\Big(\frac{9}{4p(p-1)}\Big)^{1/4}\,n^{-1/4} V^n_{pnt}\right)_{0\leq t\leq 1}
\build{\la}_{n\to\infty}^{\rm(d)} (\ov{\bf e}_t,\ov Z_t)_{0\leq t\leq 1}.
\ee
in the sense of weak convergence of the laws in the space of probability
measures on $C([0,1],\R^2)$.
\end{theorem}

The case $p=2$ of Theorem \ref{invar1} is a special case of Theorem 2.1
in \cite{LGinvar}, which is itself a
conditional version of invariance principles relating 
discrete snakes to the Brownian snake \cite{JM}. See the discussion in Section 8 of \cite{LGinvar}. 
Similar results were obtained before
by Chassaing and Schaeffer \cite{CS}. In the general case, Theorem \ref{invar1}
is a consequence of Theorem 3.3 in \cite{We}, and is also
closely related to Theorem 11 in \cite{MaMi}. 

Although Theorem \ref{invar1} will be our main tool, we will also need
another asymptotic result, which does not easily follow from Theorem \ref{invar1}
but fortunately can be deduced from the results in \cite{MaMi}. Let
$M_n$ be the random element of ${\cal M}^p_n$ that corresponds to
$\theta_n$ via the Bouttier-Di Francesco-Guitter bijection. Obviously $M_n$
is uniformly distributed over ${\cal M}^p_n$. 
Conditionally on $M_n$, let us choose
a vertex $Y_n$ of $M_n$ uniformly at random.
The pair $(M_n,Y_n)$ is then uniformly
distributed over the set of all rooted and pointed $2p$-angulations
with $n$ faces. Theorem 3 (iii) of \cite{MaMi} gives
precise information
about the profile of distances to the point $Y_n$ in the
map $M_n$ (to be precise, \cite{MaMi} imposes a special constraint on the orientation
of the root edge depending on the distinguished point in the map, but since every rooted and pointed map with
this constraint corresponds exactly to two unconstrained
rooted and pointed maps, the results of \cite{MaMi} immediately carry over to
our setting). In our special
situation, we can restate this result as follows. We write $d_{n}$ for the graph distance on the set ${\bf
m}_n$ of vertices of
$M_n$, and for every $R>0$ and $x\in {\bf m}_n$ we denote by $B_n(x,R)$
the closed ball with radius $R$ centered at $x$ in the metric space
$({\bf m}_n,d_n)$. 

\begin{proposition} 
\label{pointedmap}
For every $\alpha,\beta> 0$,
$$P\Big[\frac{1}{(p-1)n} \#B_n(Y_n,\alpha n^{1/4})
\geq \beta \Big]
\build{\la}_{n\to\infty}^{} P\Big[\ov{\cal I}\Big(
\Big[0,\Big(\frac{9}{4p(p-1)}\Big)^{1/4}\alpha\Big]\Big)\geq \beta\Big].
$$
\end{proposition}

Since $Y_n$ is uniformly distributed over ${\bf m}_n$ and $\#({\bf m}_n)=(p-1)n+2$, the convergence of
the proposition can be restated as follows. For every $\alpha,\beta>0$,
\begin{equation}
\label{pointed}
E\Big[\frac{1}{(p-1)n}\#\{y\in{\bf m}_n: \#B_n(y,\alpha n^{1/4})
\geq \beta\,(p-1)n\} \Big]
\build{\la}_{n\to\infty}^{} P\Big[\ov{\cal I}\Big(
\Big[0,\Big(\frac{9}{4p(p-1)}\Big)^{1/4}\alpha\Big]\Big)\geq \beta\Big].
\end{equation}

\section{Main results}

Recall the notation introduced in the previous section. In particular,
$M_n$ is a random rooted $2p$-angulation which is uniformly distributed over the set
$\m^p_n$, $\mm_n$ denotes the
set of vertices of $M_n$, and $\theta_n=(\tau_n,(\ell^n_u)_{u\in\tau^\circ_n})$ is the
random mobile corresponding to $M_n$ via the Bouttier-Di Francesco-Guitter bijection.
We constantly use the identification
$$\mm_n=\tau^\circ_n \cup\{\partial_n\}$$
where $\partial_n$ is the root vertex of $M_n$. The graph distance on $\mm_n$
is denoted by $d_n$. In particular, if $a,b\in\tau^\circ_n$, $d_n(a,b)$ denotes the
graph distance between $a$ and $b$ viewed as vertices in
the map $M_n$.

As in subsection 2.5, $C^n$ and $V^n$ are respectively the contour function
of the tree $\tau^\circ_n$ and the spatial
contour function of $\theta_n$. 

Following subsection 2.2, the equivalence relation $\sim_n$ on $[pn]=\{0,1,\ldots,pn\}$
is defined by declaring that $i \sim_n j$ if and only if
the $i$-th vertex in the search-depth sequence of $\tau^\circ_n$ 
is the same as the $j$-th vertex in the same sequence. Recall that this implies
$$C^n_i=C^n_j=\inf_{i\wedge j\leq k\leq i\vee j} C^n_k .$$
The quotient set $[pn]\,/\!\sim_n$ is then canonically identified
with $\tau^\circ_n$ and thus with the set of vertices of
$M_n$ other than the root $\partial_n$. If $a\in \tau^\circ_n$ and
$i\in [pn]$, we will abuse notation by writing $a\sim_n i$
if $i$ is a representative of $a$ viewed as an element of
$[pn]\,/\!\sim_n$ (similar abuses of notation will occur
for other equivalence relations). With this notation,
if $a\sim_n i$, we have $d_n(\partial_n,a)=\ell^n_a=V^n_i$, by the properties
the Bouttier-Di Francesco-Guitter bijection. If $i,j\in [pn]$ and $a,b\in \tau^\circ_n$
are such that $a\sim_n i$ and $b\sim_n j$, we will also
write $d_n(i,j)=d_n(a,b)$. 

For every $i,j\in[pn]$, we put
$$d_n^\circ(i,j)=
V^n_i+V^n_j -2\inf_{i\wedge i\leq k\leq i\vee j} V^n_k +2.$$

\begin{lemma}
\label{simple-bound}
For every $i,j\in [pn]$,
$$d_n(i,j)\leq d_n^\circ(i,j).
$$
\end{lemma}

\proof Fix $i\in [pn]$ and let $a\in \tau^\circ_n$ be such that
$a\sim_n i$. Let $q=V^n_i=d_n(\partial_n,a)$. We set $i_q=i$
and for every $k\in\{1,\ldots,q-1\}$,
$$i_k=\inf\{\ell \geq i:V^n_\ell = k\}.$$
 From the construction of edges in the Bouttier-Di Francesco-Guitter bijection,
it is immediate to see that $d_n(i_k,i_{k-1})=1$
for every $2\leq k\leq q$.

We also fix $j\in[pn]$ and let $b\in \tau^\circ_n$ be such that
$b\sim_n j$, and we set $r=V^n_j=d_n(\partial_n,b)$. We define
similarly the sequence $j_r=j,j_{r-1},\ldots,j_1$. Then:
\begin{description}
\item{$\bullet$} Either $\inf_{i\wedge i\leq k\leq i\vee j} V^n_k=1$
and the bound of the lemma is just the triangle inequality
$d_n(i,j)=d_n(a,b)\leq d_n(\partial_n,a)+d_n(\partial_n,b)$.
\item{$\bullet$} Or $\inf_{i\wedge i\leq k\leq i\vee j} V^n_k=\ell\geq 2$,
and we have $i_{\ell-1}=j_{\ell-1}$. The bound of the lemma
follows by writing:
$$d_n(i,j)\leq d_n(i_{\ell-1},i_q)+d_n(j_{\ell-1},j_r)\leq q+r-2\ell +2.$$
\end{description}
\par\cq

We extend the definition of
$d_n(i,j)$ and $d_n^\circ(i,j)$ to noninteger values
of $i$ and $j$ by linear interpolation. If $s,t\in[0,pn]$,
we set
\begin{eqnarray*}
d_n(s,t)&=&(s-\lfloor s\rfloor)(t-\lfloor
t\rfloor)d_n(\lceil s\rceil,
\lceil t\rceil) +(s-\lfloor s\rfloor)(\lceil t\rceil -t)d_n(\lceil
s\rceil,
\lfloor t\rfloor)\cr
&+& (\lceil s\rceil-s)(t-\lfloor
t\rfloor)d_n(\lfloor s\rfloor,
\lceil t\rceil)
+(\lceil s \rceil -s)(\lceil t\rceil -t)d_n(\lfloor
s\rfloor,
\lfloor t\rfloor),
\end{eqnarray*}
with the notation $\lfloor s\rfloor=\sup\{k\in\Z:k\leq s\}$
and $\lceil s\rceil=\inf\{k\in\Z:k>s\}$. We 
define $d_n^\circ(s,t)$ in a similar way. Obviously the bound
$d_n(s,t)\leq d_n^\circ(s,t)
$ remains valid for reals $s,t\in[0,pn]$. Furthermore, the triangle inequality
$d_n(s,u)\leq d_n(s,t)+d_n(t,u)$ also holds for every
$s,t,u\in [0,pn]$.

As a straightforward consequence of (\ref{basicinvar}) and the definition
of $d^\circ_n(s,t)$, we have
\be
\label{basic-dist}
\left(
\Big(\frac{9}{4p(p-1)}\Big)^{1/4}\,n^{-1/4}
d^\circ_n(pns,pnt)\right)_{0\leq s\leq 1, 0\leq t\leq 1}
\build{\la}_{n\to\infty}^{\rm(d)} (D^\circ(s,t))_{0\leq
s\leq 1, 0\leq t\leq 1}
\ee
where
$$D^\circ(s,t)
=\ov Z_s+ \ov Z_t -2 \inf_{s\wedge t\leq r \leq s\vee t} \ov Z_r$$
and the limit holds in the sense of weak convergence
in the space of probability measures on $C([0,1]^2,\R)$.

\begin{proposition}
\label{tightness}
The sequence of the laws of the processes
$$\left(n^{-1/4}\,d_n(pns,pnt)\right)_{0\leq s\leq 1,0\leq t\leq 1}$$
is tight in the space of probability measures on
$C([0,1]^2,\R)$. Let $\bf C$ be the space of isometry classes of
compact metric spaces, which is equipped with the Gromov-Hausdorff measure.
The sequence of the laws of the metric spaces $(\mm_n,n^{-1/4}d_n)$
is tight in the space of probability measures on $\bf C$.
\end{proposition}

\proof First observe that, for every $s,t,s',t'\in[0,1]$,
\begin{eqnarray}
\label{tight0}
|n^{-1/4}\,d_n(pns,pnt)-n^{-1/4}\,d_n(pns',pnt')|
&\leq&n^{-1/4}(d_n(pns,pns')+d_n(pnt,pnt'))\nonumber\\
&\leq&n^{-1/4}(d^\circ_n(pns,pns')+d^\circ_n(pnt,pnt')).
\end{eqnarray}
 From the convergence (\ref{basic-dist}), we have for every
$\delta,\varepsilon>0$,
\begin{equation}
\label{tight1}
\limsup_{n\to \infty} P\left(\sup_{|s-s'|\leq \delta}
n^{-1/4}\,d^\circ_n(pns,pns') \geq \varepsilon\right)
\leq P\left(\sup_{|s-s'|\leq \delta} D^\circ(s,s')\geq
\Big(\frac{4p(p-1)}{9}\Big)^{1/4}\,\varepsilon\right).
\end{equation}
Let $\eta>0$ and for every $k\geq 1$ set $\varepsilon_k=2^{-k}$.
We apply (\ref{tight1}) with $\varepsilon=\varepsilon_k$ and note that
we can then choose $\delta_k>0$ sufficiently small so that the
right-hand side of (\ref{tight1}) is strictly less than $2^{-k}\eta$.
Therefore, there exists an integer $n_k$ such that,
for every $n\geq n_k$,
\be
\label{tight2}
P\left(\sup_{|s-s'|\leq \delta_k}
n^{-1/4}\,d^\circ_n(pns,pns') \geq \varepsilon_k\right)
\leq 2^{-k}\eta.
\ee
By choosing $\delta_k$ even smaller if necessary, we may assume that
(\ref{tight2}) holds for every $n\geq 1$. It follows that,
for every $n\geq 1$,
\be
\label{tight3}
P\left(\bigcap_{k\geq 1}\left\{\sup_{|s-s'|\leq \delta_k}
n^{-1/4}\,d^\circ_n(pns,pns') \leq \varepsilon_k\right\}\right)
\geq 1-\eta.
\ee
Let $K$ denote the set of all functions $\omega\in C([0,1]^2,\R)$
such that $\omega(0,0)=0$ and, for every $k\geq 1$,
$$\sup\{|\omega(s,t)-\omega(s',t')|:|s-s'|\leq \delta_k,|t-t'|\leq
\delta_k\}\leq 2\,\varepsilon_k.$$
Then $K$ is a compact subset of $C([0,1]^2,\R)$. By (\ref{tight0})
and (\ref{tight3}), the probability that the random function
$(s,t)\la n^{-1/4}d_n(pns,pnt)$ belongs to $K$ is bounded
below by $1-\eta$, for every $n\geq 1$. Since $\eta$
was arbitrary, this completes the proof of the first assertion. 

The second assertion is an easy consequence of the first one and
the Gromov compactness
criterion (Theorem 7.4.15 in \cite{BBI}). We omit details, since this
result is not really needed in what follows.
\cq

\medskip
 From (\ref{basicinvar}) and Proposition \ref{tightness}, there exists a
strictly increasing sequence
$(n_k)_{k\geq 1}$ such that along this sequence we have
the joint convergence in distribution
\begin{eqnarray}
\label{basic}
&&\hspace{-6mm}\left({1\over 2}\sqrt{\frac{p}{p-1}}\,n^{-1/2}\,C^n_{pnt},
\Big(\frac{9}{4p(p-1)}\Big)^{1/4}\,n^{-1/4} V^n_{pnt},
\Big(\frac{9}{4p(p-1)}\Big)^{1/4}\,n^{-1/4}\,d_n(pns,pnt)\!\right)
_{0\leq s\leq 1, 0\leq
t\leq 1}\nonumber\\
\noalign{\smallskip}
&&\qquad\build{\la}_{n\to\infty}^{} \left(\ov{\bf e}_t,\ov Z_t
,D(s,t)\right)_{0\leq s\leq 1,0\leq t\leq
1}.
\end{eqnarray}
Here the limiting triple  $\left(\ov{\bf e}_t,\ov Z_t
,D(s,t)\right)$ is defined on a suitable probability space, the pair
$(\ov{\bf e},\ov Z)$ obviously has the same distribution as before, and
$D$ is a continuous process indexed by $[0,1]^2$ and taking values
in $\R_+$. In the remaining part of this work, we restrict our attention
to values of $n$ belonging to the sequence $(n_k)_{k\geq 1}$.
In particular, when we pass to the limit as $n\to \infty$, this
always means along the sequence $(n_k)_{k\geq 1}$.

Thanks to the Skorokhod representation theorem, we may and will assume
that the convergence (\ref{basic}) holds almost surely, in the
sense of uniform convergence over $[0,1]^2$.
Strictly speaking, we should replace for every $n\geq 1$ the
random mobile $\theta_n$ (respectively the random map $M_n$)
with another random mobile $\wt\theta_n$
(resp. another random map $\wt M_n$)
having the same distribution, but we do not keep
track of this replacement in the notation.

The next proposition records some properties of the
random function $D(s,t)$. We write $\simeq$ instead
of $\simeq_{\ov\eg}$
for the equivalence relation defining the tree
$\t_{\ov{\bf e}}\;$: $\t_{\ov{\bf e}}=[0,1]\,/\!\simeq$
as was explained in subsection 2.3.

\begin{proposition}
\label{prop-distance} The following properties hold almost surely.
\par\noindent{\rm(i)} For every $s,t,u\in[0,1]$,
\ba
&&D(s,s)=0\\
&&D(s,t)=D(t,s)
\ea
and
$$D(s,u)\leq D(s,t) + D(t,u).$$
{\rm (ii)} For every $s,t\in [0,1]$,
$$D(s,t)\leq D^\circ(s,t).$$
{\rm(iii)} For every $s,t\in [0,1]$, the property
$s \simeq t$ implies $D(s,t)=0$.
\par\noindent{\rm (iv)} For every $s\in[0,1]$, $D(0,s)=\ov Z_s$.
\end{proposition}

\proof Except for the first one, the properties in (i) are immediate from
the analogous properties for $d_n$ and the (almost sure) convergence
(\ref{basic}). Similarly, (ii) follows from Lemma \ref{simple-bound} and
the convergence (\ref{basic-dist}), which holds a.s. along the
sequence $(n_k)_{k\geq 1}$ if (\ref{basic}) also holds a.s.
along this sequence. The first
property in (i) then readily follows from (ii).

Let us prove (iii). Let $s,t\in[0,1]$ with $s<t$. If
$s\simeq t$, we have
$$\ov{\bf e}_s=\ov{\bf e}_t=\inf_{s\leq r\leq t} \ov{\bf e}_r.$$
Suppose first that $\ov{\bf e}_r> \ov{\bf e}_s$ for every
$r\in]s,t[$. From the uniform convergence of the function
${1\over 2}\sqrt{{p}/{(p-1)}}\,n^{-1/2}\,C^n_{pnt}$ towards
$\ov{\bf e}_t$, an elementary argument yields the existence of
two sequences $(i_n)$ and $(j_n)$ of integers in $[pn]$ such that:
\begin{description}
\item{$\bullet$} $\displaystyle\frac{i_n}{pn}\la s$ and $
\displaystyle\frac{j_n}{pn}\la
t$ as $n\to \infty$.
\item{$\bullet$} For $n$ sufficiently large, $j_n\geq i_n+ 2$
and $C^n_{i_n}=C^n_{j_n}<{\displaystyle\inf_{i_n<k<j_n} C^n_k}$.
\end{description}
As we already noticed in subsection 2.2, the last property ensures that $i_n\sim_n j_n$ and thus
$d_n(i_n,j_n)=0$. By passing to the limit $n\to\infty$, we get
$D(s,t)=0$.

If $\ov{\bf e}_r= \ov{\bf e}_s$ for some
$r\in]s,t[$, then $r$ is necessarily unique, because otherwise
the tree $\t_{\ov{\bf e}}$, which is isometric to $\t_\eg$, would have a point with multiplicity
strictly greater than $3$. By the preceding argument,
$D(s,r)=D(r,t)=0$ and thus $D(s,t)=0$ by the triangle inequality
in (i).

Let us finally prove (iv). Let $s\in[0,1]$ and let $(i_n)$ be
a sequence of integers such that $i_n/(pn) \la s$ as $n\to\infty$.
From the properties of the
Bouttier-Di Francesco-Guitter bijection, we know that 
$d_n(0,i_n)=V^n_{i_n}$. On the other hand, 
(\ref{basic}) ensures that $(9/(4p(p-1)))^{1/4}n^{-1/4}d_n(0,i_n)$
converges to $D(0,s)$, and that $(9/(4p(p-1)))^{1/4}n^{-1/4}V^n_{i_n}$
converges to $\ov Z_s$. The desired result follows.
\cq

\smallskip
We define an equivalence relation
$\approx$ on $[0,1]$ by setting
$$s\approx t\quad\hbox{if and only if}\quad D(s,t)=0.$$
Clearly, $D$ induces a metric, which we still denote by $D$, on the
quotient set $[0,1]\,/\!\approx$. The bound $D\leq D^\circ$ ensures that the canonical projection
from $[0,1]$ onto $[0,1]\,/\!
\approx$ is continuous when
$[0,1]\,/\!
\approx$ is equipped with the metric $D$. In particular the metric space
$([0,1]\,/\!\approx, D)$ is compact.

For our purposes, it will be convenient to view this metric space as
a quotient of the real tree $\t_{\ov{\bf e}}$.
By property (iii) of the previous proposition, we may define
$D(a,b)$ for $a,b\in \t_{\ov{\bf e}}=[0,1]\,/\!\simeq$ simply by setting
$D(a,b)=D(s,t)$ where $s$, resp. $t$, is any representative of $a$, resp.
$b$, in $[0,1]$. The equivalence relation $\approx$
then makes sense on $\t_{\ov{\bf e}}$, and the quotient space
$(\t_{\ov{\bf e}}\,/\!
\approx,D)$ is obviously isometric to $([0,1]\,/\!\approx, D)$.
As a consequence of Proposition \ref{prop-distance} (iv)
and the triangle inequality, 
for every $a,b\in\t_{\ov{\bf e}}$, the condition $D(a,b)=0$
implies $\ov Z_a=\ov Z_b$.

Before stating the main result, we need to introduce some
additional notation. For every $a,b\in \t_{\ov{\bf e}}$, we set
$$D^\circ(a,b)=\inf\{D^\circ(s,t):s,t\in[0,1], a\simeq s, b\simeq t\}.$$
Suppose that neither of the relations $a\prec b$
and $b\prec a$ holds, and assume for definiteness
that $a<b$. Then the infimum in the definition of $D^\circ(a,b)$
is attained when $[s,t]$ is the minimal subinterval of $[0,1]$ such
that $a\simeq s$ and $b\simeq t$, and it follows that
$$D^\circ(a,b)=\ov Z_a +\ov Z_b - 2 \inf_{c\in [a,b]} \ov Z_c$$
where $[a,b]$ is the lexicographical interval between $a$
and $b$ in $\t_{\ov {\bf e}}$, as defined in subsection 2.3. On the other hand, if
$a\prec b$, then the preceding formula does not necessarily hold: We
have instead
$$D^\circ(a,b)=\ov Z_a +\ov Z_b - 2 \sup\left(\inf_{c\in [a,b]} \ov Z_c,
\inf_{c\in [b,a]} \ov Z_c\right).$$
The function $D^\circ(a,b)$, $a,b\in\t_{\ov{\bf e}}$ needs not satisfy the
triangle inequality. For this reason, we set
for every  $a,b\in\t_{\ov{\bf e}}$,
$$D^*(a,b)=\inf\left\{\sum_{i=1}^q D^\circ(a_{i-1},a_i)\right\}$$
where the infimum is over all choices of the integer $q\geq 1$
and of the finite sequence $a_0,a_1,\ldots,a_q$
in $\t_{\ov{\bf e}}$ such that $a_0=a$ and $a_q=b$.

Since $D\leq D^\circ$, and $D$ satisfies the triangle inequality, it
is clear that we have
$$0\leq D(a,b)\leq D^*(a,b)\leq D^\circ (a,b)$$
for every $a,b\in\t_{\ov{\bf e}}$.

We can now state our main result. Recall that we are restricting
our attention to values of $n$ belonging to the sequence
$(n_k)_{k\geq 1}$, and that we assume that the convergence
(\ref{basic}) holds a.s. along this sequence.

\begin{theorem}
\label{main}
We have almost surely
$$\left({\bf m}_n,\Big(\frac{9}{4p(p-1)}\Big)^{1/4}\,n^{-1/4} d_n\right)
\build{\la}_{n\to\infty}^{} (\t_{\ov{\bf e}}\,/\!\approx,D)$$
in the sense of the Gromov-Hausdorff distance on compact
metric spaces. In addition, a.s. for every $a,b\in\t_{\ov{\bf e}}$,
the relation $a\approx b$ holds if and only if one of the
following equivalent properties holds:
\begin{description}
\item{\rm(i)} $D(a,b)=0$.
\item{\rm(ii)} $D^*(a,b)=0$.
\item{\rm(iii)} $D^\circ(a,b)=0$.
\end{description}
\end{theorem}

\noindent{\bf Remarks.} (a)
Although the process $D$ may depend on the sequence $(n_k)_{k\geq 1}$,
the equivalence relation $\approx$ does not, since it can be defined
by either (ii) or (iii) in Theorem \ref{main}.
As was already observed in the introduction, this guarantees that the
limiting compact metric space $(\t_{\ov{\bf e}}\,/\!\!\approx,D)$
is homeomorphic to $(\t_{\ov{\bf e}}\,/\!\approx,D^*)$, and thus that
its topology does not depend on the choice of the
sequence $(n_k)_{k\geq 1}$ (nor on the value of $p$). Still it is tempting
to conjecture that $D(a,b)=D^*(a,b)$, for every $a,b\in\t_{\ov\eg}$.
If this conjecture is correct, the convergence (\ref{basic}), or
that of Theorem \ref{main}, does not require the use of
a subsequence.

\noindent(b) It is not hard to prove that equivalence classes in $\t_{\ov\eg}$ for the 
equivalence relation $\approx$ can contain only $1$, $2$ or $3$ points. For
every fixed $s\in[0,1]$,
it is easy to verify that the equivalence 
class of $p_{\ov\eg}(s)$ is a singleton a.s. Furthermore,
one can check that a.s. for every rational numbers $r,s,t,u$
such that $0\leq r<s<t<u\leq 1$ one has 
$$\inf_{r\leq x\leq s} \ov Z_x \ne \inf_{t\leq x\leq u} \ov Z_x.$$
(The easiest way to derive this property is to use the Brownian snake approach
that is presented below in Section 5.) It follows that an equivalence class
cannot contain more than $3$ points. Conversely, if
we are given two rationals $0<r<s<1$, there exists an a.s. unique $y\in]r,s[$
such that
$$\ov Z_y=\inf_{r\leq x\leq s} \ov Z_x,$$
and the vertex of $\t_{\ov\eg}$ corresponding to $y$ is a leaf of $\t_{\ov\eg}$.
Set $t_1=\sup\{u\leq r:\ov Z_u=\ov Z_y\}$ and $t_2=\inf\{u\geq s:\ov Z_u=\ov Z_y\}$.
Then $t_1\approx y\approx t_2$, and $t_1,y$ and $t_2$
correspond to different vertices of the tree $\t_{\ov{\bf e}}$. To summarize, the equivalence class of a
typical vertex $a\in\t_{\ov\eg}$ is a singleton, but there is a continuum of equivalence
classes consisting of pairs, and there are countably many equivalence classes
containing three elements. These properties are not used below. They will be derived
in greater detail
in the subsequent paper \cite{LGP} where they play an important role.

\medskip
\noindent{\bf Proof of Theorem \ref{main} (first part):} The main difficulty in the proof of Theorem
\ref{main} comes from the implication (i)$\Rightarrow$(iii). Notice that the other
implications (iii)$\Rightarrow$(ii)$\Rightarrow$(i) are trivial. The
implication (i)$\Rightarrow$(iii) is established in the next section.
We now prove the first assertion of
Theorem \ref{main}.

Recall that the metric spaces
$(\t_{\ov{\bf e}}\,/\!\approx,D)$ and $([0,1]\,/\!\approx,D)$
are isometric. For every integer $n$, consider the equivalence
relation $\approx_n$ defined on $[0,1]$ by setting
$$s\approx_n t\hbox{\quad if and only if\quad} d_n(\lfloor pns\rfloor,
\lfloor pnt\rfloor)=0.$$
Clearly, the quotient space $E_n:=[0,1]\,/\!\!\approx_n$
equipped with the metric $\delta_n(s,t)=d_n(\lfloor pns\rfloor,
\lfloor pnt\rfloor)$ is isometric to $(\tau_n^\circ,d_n)$
or equivalently to $({\bf m}_n\backslash\{\partial_n\},d_n)$.

Since
$$d_{GH}(({\bf m}_n\backslash\{\partial_n\},n^{-1/4}d_n),({\bf
m}_n,n^{-1/4}d_n))
\build{\la}_{n\to\infty}^{} 0$$
the first part of Theorem \ref{main} reduces to checking that
we have a.s.
\be
\label{convDGH}
d_{GH}\left(\Big(E_n,
\Big(\frac{9}{4p(p-1)}\Big)^{1/4}\,n^{-1/4}\delta_n\Big),(E_\infty,D)\right)
\build{\la}_{n\to\infty}^{} 0
\ee
where $E_\infty:=[0,1]\,/\!\approx$.

To this end, we construct a correspondence between the
metric spaces $E_n$ and $E_\infty$ by setting
$${\mathcal C}_n=\{(a,b)\in E_n\times E_\infty:
\hbox{there exists }t\in[0,1]\hbox{ such that } a\approx_n t
\hbox{ and }b\approx t\}.$$
In order to bound the distortion of this correspondence,
consider two pairs $(a,b),(a',b')\in{\mathcal C}_n$. By definition,
there exist $s,t\in[0,1]$
such that $a\approx_n s,b\approx s$ and $a'\approx_n t,
b'\approx t$. Then we have
\ba
&&\delta_n(a,a')=d_n(\lfloor pns\rfloor,\lfloor pnt\rfloor)\\
&&D(b,b')=D(s,t).
\ea
Thus, when $E_n$ is equipped with the distance
$({9}/{(4p(p-1))})^{1/4}\,n^{-1/4}\delta_n$, and $E_\infty$ with the
distance $D$, the distortion of ${\mathcal C}_n$ is
\ba
&&\sup_{(a,b),(a',b')\in{\mathcal C}_n}
\Big|\Big(\frac{9}{4p(p-1)}\Big)^{1/4}\,n^{-1/4}\delta_n(a,a')
-D(b,b')\Big|\\
&&\qquad\leq \sup_{s,t\in[0,1]}
\Big|\Big(\frac{9}{4p(p-1)}\Big)^{1/4}\,n^{-1/4}
d_n(\lfloor pns\rfloor,\lfloor pnt\rfloor)-
D(s,t)\Big|,
\ea
which tends to $0$ a.s. by (\ref{basic}). The first assertion of Theorem \ref{main} now follows
from the known result connecting the Gromov-Hausdorff
distance between two compact metric spaces with the infimum
of the distortion of correspondences between these two spaces (Theorem 7.3.25
in \cite{BBI}).
\cq 

\medskip
Before proceeding to the second part of the proof of Theorem \ref{main}, let
us state and prove a closely related result.

\begin{proposition}
\label{finite-margi}
Let $k\geq 1$ be an integer. For every $n\geq 1$, let $Y^n_1,\ldots,Y^n_k$ be 
$k$ random variables which conditionally given $M_n$ are independent and
uniformly distributed over ${\bf m}_n$. Also, given the triple $(\ov\eg,\ov Z,D)$, let
$Y^\infty_1,\ldots,Y^\infty_k$ be random variables
with values in $\t_{\ov\eg}$ which are independent and distributed according to $\lambda$.
Then,
$$\left(\Big(\frac{9}{4p(p-1)}\Big)^{1/4}\,n^{-1/4}
\,d_n(Y^n_i,Y^n_j)\right)_{1\leq i\leq k,1\leq j\leq k}
\build{\la}_{n\to\infty}^{\rm (d)}
\left(D(Y^\infty_i,Y^\infty_j)\right)_{1\leq i\leq k,1\leq j\leq k}.$$
\end{proposition}

\noindent{\bf Remarks.} (a) Informally, Proposition \ref{finite-margi} means that the convergence in
Theorem \ref{main} can be reinforced in the sense of convergence of measured metric spaces, provided ${\bf
m}_n$ is equipped with the uniform probability measure and $\t_{\ov{\bf e}}\,/\!\approx$
is equipped with the image of $\lambda$ under the canonical projection. We could give other 
versions of this reinforcement: See Chapter $3\frac{1}{2}$ of the book \cite{Gro}
for various notions of convergence of measured metric spaces. Here we content ourselves
with the preceding proposition, which will be useful in Section 6 below.

\noindent (b) The reader may be puzzled by our 
assumption on $Y^\infty_1,\ldots,Y^\infty_k$, since we seem
to be dealing with random variables taking values in a {\it random} state space. It is however 
a straightforward matter to give a mathematically rigorous (although less intuitive)
version of the statement of the proposition.

\medskip
\proof Recall that ${\bf m}_n=\tau^\circ_n \cup\{\partial_n\}$. We may
and will assume that $Y^n_1,\ldots,Y^n_k$ are 
uniformly distributed over $\tau^\circ_n$ rather than over ${\bf m}_n$.

Then let $U_1,\ldots,U_k$ be $k$ independent random variables which are uniformly
distributed over $[0,1]$ and independent of all other random quantities we have
considered until now. We may then take $Y^\infty_i=p_{\ov\eg}(U_i)$
for every $1\leq i\leq k$. Also, for every $1\leq i\leq k$, we let
$\wt Y^n_i$ be the equivalent class of $\lfloor pnU_i \rfloor$ in the quotient set
$[pn]\,/\!\sim_n=\tau^\circ_n$. The (almost sure) convergence (\ref{basic}) implies that
$$\left(\Big(\frac{9}{4p(p-1)}\Big)^{1/4}\,n^{-1/4}
\,d_n(\wt Y^n_i,\wt Y^n_j)\right)_{1\leq i\leq k,1\leq j\leq k}
\build{\la}_{n\to\infty}^{\rm (a.s.)}
\left(D(Y^\infty_i,Y^\infty_j)\right)_{1\leq i\leq k,1\leq j\leq k}.$$
This does not immediately give us the desired result, because the 
variables $\wt Y^n_i$ are not uniformly distributed over $\tau^\circ_n$.
Still we will see that in a sense they are close enough to variables that have 
the desired uniform distribution. To this end, for every $n$ and every $1\leq i\leq k$, set
$$k^n_i=\lceil ((p-1)n+2)U_i\rceil$$
and let $Y^n_i$ be the $k^n_i$-th element in the sequence of vertices of $\tau_n^\circ$
listed in lexicographical order. Clearly, the variables $Y^n_i$ have the properties
stated in the proposition. To complete the proof, it is therefore enough to
check that, for every $1\leq i\leq k$,
\be
\label{margitech1}
n^{-1/4}d_n(Y^n_i,\wt Y^n_i)\build{\la}_{n\to\infty}^{\rm a.s.} 0.
\ee
Note that a.s. for every $t\in]0,1]$, the number of distinct vertices of $\tau^\circ_n$
that appear in the search-depth sequence before rank $\lfloor pnt\rfloor$ behaves
as $(p-1)nt$ when $n\to\infty$. To see this, observe 
 that in the evolution of the contour function
of $\tau^\circ_n$ each step which is not downwards corresponds in the 
search-depth sequence to a vertex of $\tau^\circ_n$ that has not been visited
before, and then use (\ref{basic}) to see that the number of downward steps
before time $\lfloor pnt\rfloor$ behaves like $nt$ when $n\to\infty$ (indeed the difference between
the numbers of upward and downward steps is $O(n^{1/2})$ as $n\to\infty$).

From the preceding remarks, we get that a.s. for every $s,t\in[0,1]$
such that $s<U_i<t$, if $n$ is large enough, the vertex $Y^n_i$
is visited by the search-depth sequence during the time interval $[\lfloor pns\rfloor,\lfloor pnt\rfloor]$.
Thus, for $n$ sufficiently large,
$$n^{-1/4}d_n(Y^n_i,\wt Y^n_i)\leq n^{-1/4}\sup_{\lfloor pns\rfloor\leq j\leq \lfloor pnt\rfloor}
d_n(j,\lfloor pnU_i\rfloor).$$
The right-hand side can be made arbitrarily small when $n$ is large by choosing $s$ and $t$
sufficiently close to $U_i$. This completes the proof of
(\ref{margitech1}) and of Proposition \ref{finite-margi}. \cq

\section{The key step}

This section is devoted to the second part of the proof of Theorem \ref{main}, that is 
to the proof of the implication
(i)$\Rightarrow$(iii) in this theorem. We start with a lemma.

\begin{lemma}
\label{linesegment}
Almost surely, for every $a,b\in\t_{\ov{\bf e}}$, the condition
$D(a,b)=0$ implies $\ov Z_c\geq \ov Z_a=\ov Z_b$
for every $c\in\llbracket a,b\rrbracket$.
\end{lemma}

\proof We already noticed that the condition $D(a,b)=0$ forces
$\ov Z_a=\ov Z_b$. We can immediately exclude the case $a=\ov\rho$
since this would imply $\ov Z_b=\ov Z_a=0$
and $b=\ov \rho=a$. Then we can assume without loss of
generality that $a<b$. We argue by contradiction, assuming that
there exists $c\in\rrbracket a,b \llbracket$ such that 
$\ov Z_c< \ov Z_a$. For definiteness, we assume that 
$c\in \llbracket a\wedge b,a\llbracket$. The symmetric case
$c\in \llbracket a\wedge b,b\llbracket$ is treated in a 
similar manner.

Let $s<t$ be such that $a\simeq s$ and $b\simeq t$. 
We can then find $r\in]s,t[$ such that $c\simeq r$. Choose
$i_n,j_n,k_n\in[pn]$,
with $i_n\leq k_n\leq j_n$, such that $i_n/(pn)\la s$,
$j_n/(pn)\la t$ and $k_n/(pn)\la r$.
Denote by $a_n,b_n,c_n$ the vertices in $\tau^\circ_n$
corresponding respectively to $i_n,j_n,k_n$. Since 
$c\in \llbracket a\wedge b,a\llbracket$, a simple argument using the
convergence of the first components in (\ref{basic}),
and the remarks of the beginning of subsection 2.2,
shows that $k_n$ can be chosen in such a way that
$c_n\in\rrbracket a_n\wedge b_n,a_n\llbracket$ for every
$n$ sufficiently large. Denote by $\tau^\circ_n(c_n)$
the set of all descendants of $c_n$ in $\tau^\circ_n$. Then
$a_n\in\tau^\circ_n(c_n)$ but $b_n\notin\tau^\circ_n(c_n)$.

By (\ref{basic}) and our assumption $D(a,b)=0$
we know that $d_n(a_n,b_n)=o(n^{1/4})$ as $n\to\infty$. Let
$\gamma_n=(\gamma_n(i),0\leq i\leq d_n(a_n,b_n))$ be a geodesic path
from $a_n$ to $b_n$ in the map $M_n$. When $n$ is large,
the path $\gamma_n$ must lie entirely in $\tau^\circ_n$,
because if $\partial_n$ belongs to this path the
equality $d_n(a_n,b_n)=d_n(a_n,\partial_n)+ d_n(\partial_n,b_n)=V^n_{i_n}+V^n_{j_n}$
yields a contradiction with the property $d_n(a_n,b_n)=o(n^{1/4})$.

Denote by $g_n$
the last point on the geodesic $\gamma_n$ that belongs 
to $\tau^\circ_n(c_n)$. Since $g_n$ is a point of the
geodesic $\gamma_n$ and $d_n(a_n,b_n)=o(n^{1/4})$, we have
$$\ell^n_{g_n}=d_n(\partial_n,g_n)=d_n(\partial_n,a_n)+o(n^{1/4})
=n^{1/4} \ov Z_a +o(n^{1/4})$$
as $n\to\infty$. On the other hand, since $k_n/(pn)\la r$
and $c\simeq r$,
$$\ell^n_{c_n}=V^n_{k_n}=n^{1/4}\ov Z_c +o(n^{1/4})$$
as $n\to\infty$. Hence, for $n$ large we 
must have $\ell^n_{g_n}>\ell^n_{c_n}$.

Using the way edges of the map $M_n$ are reconstructed from the 
mobile $\theta_n$, we now see that any edge starting from $g_n$
in $M_n$ connects $g_n$ with another point of $\tau^\circ_n(c_n)$.
Indeed, any successor of the vertex $g_n$ must clearly lie 
in $\tau^\circ_n(c_n)$ because in the search-depth sequence 
of $\tau^\circ_n$, a vertex with label $\ell^n_{g_n}-1$ will be
visited after the last visit of $g_n$ before coming back to $c_n$ and exiting the tree $\tau^\circ_n(c_n)$.
Similarly, $g_n$ cannot be a successor of a
vertex $h\notin\tau^\circ_n(c_n)$: If this were the case we would have
$\ell^n_h-1=\ell^n_{g_n}>\ell^n_{c_n}$, and the search-depth sequence 
of $\tau^\circ_n$ would visit a vertex with label $\ell^n_h-1$ after visiting $h$
before entering the set $\tau^\circ_n(c_n)$. Finally, the fact that
$g_n$ is not connected to any point outside $\tau^\circ_n(c_n)$
gives a contradiction with our choice of $g_n$. \cq

\begin{proposition}
\label{ancestorline}
Almost surely, for every pair $(a,b)$ in $\tree$ such that
$a$ is an ancestor of $b$ and $a\ne b$, we have $D(a,b)>0$. 
\end{proposition}

\proof We argue by contradiction, assuming that there exists 
a pair $(a,b)$ in $\tree$ such that
$a$ is an ancestor of $b$, $a\ne b$ and $D(a,b)=0$.
Notice that  the case $a=\ov\rho$ is excluded since we already
know from Proposition \ref{prop-distance} (iv) that $D(\ov\rho,b)=\ov Z_b>0$
for every $b\ne
\ov\rho$. So we assume that  $a\ne \ov\rho$.
Recall that we have automatically $\ov Z_a=\ov Z_b$.

Let $s,t\in[0,1]$ be such that 
$a\simeq s$ and $b\simeq t$. Since $a$ is an ancestor of $b$
we can choose $s$ and $t$ such that $s<t$ and $\ov{\bf e}_r>\ov{\bf e}_s$
for every $r\in]s,t]$.
Since
$D(s,t)=D(a,b)=0$,
(\ref{basic}) gives
$$n^{-1/4}d_n(pn s,pnt)\build{\la}_{n\to\infty}^{} 0.$$
So, for every $n$, we can find $i^\circ_n,j_n\in[pn]$ such that
$i^\circ_n\leq j_n$,
$|i^\circ_n-pns|\leq 1$, $|j_n-pnt|\leq 1$ and
$$n^{-1/4}d_n(i^\circ_n,j_n)\build{\la}_{n\to\infty}^{} 0.$$
Let $i_n=\sup\{k\in[i^\circ_n,j_n]\cap\Z: C^n_k=C^n_{i^\circ_n}\}$. By 
(\ref{basic}) and the condition $\ov{\bf e}_r>\ov{\bf e}_s$
for every $r\in]s,t]$, we must have $n^{-1}(i_n-i^\circ_n)\la 0$
as $n\to \infty$. As a consequence,  $n^{-1/4}d_n(i^\circ_n,i_n)\la 0$.
Let $a_n$ and $b_n$ be the vertices in $\tau^\circ_n$ such that
$a_n\sim_n i_n$ and $b_n\sim_n j_n$. Then provided $n$ is
sufficiently large, the remarks of subsection 2.2 show that $a_n$ is an ancestor of $b_n$. Moreover
we have
\be
\label{distzero}
n^{-1/4}d_n(a_n,b_n)=n^{-1/4}d_n(i_n,j_n)\build{\la}_{n\to\infty}^{} 0.
\ee

By Lemma \ref{linesegment}, we also know that $\ov Z_c\geq \ov Z_a$
for every $c\in\llbracket a,b\rrbracket$. Recall that the conditioned
tree $(\tree, \ov Z)$ is obtained by re-rooting the unconditioned tree
$(\t_{\bf e},Z)$ at the vertex corresponding to the minimal spatial 
position, and that along a given line segment of $\t_{\bf e}$,
$Z$ evolves like linear Brownian motion. Since local minima 
of linear Brownian motion are distinct, a simple argument then 
shows that the equality $\ov Z_c=\ov Z_a$ can hold for at most
one value of $c\in\rrbracket a,b\llbracket$. Hence, we can find 
$\eta>0$ such that the properties $c\in\rrbracket a,b\llbracket$ and
$d_{\ov{\bf e}}(a,c)<\eta$ imply $\ov Z_c>\ov Z_a$.

Since $\ov{\bf e}_r>\ov{\bf e}_s$
for every $r\in]s,t]$, Lemma \ref{increasepoint} 
implies that for every $\varepsilon>0$,
$$\inf_{r\in[s,s+\varepsilon]}\ov Z_r<\ov Z_s.$$
It follows that there exists one (in fact infinitely many) subtree
$\t^1$ from the left side of $\llbracket a,b\rrbracket$, with
root $\rho^1\in\rrbracket a,b \llbracket$, such that
$d_{\ov{\bf e}}(a,\rho^1)<\eta$  and
$$\inf_{c\in\t^1} \ov Z_c <\ov Z_a.$$
We denote by $[\alpha,\beta]$ the interval coding $\t^1$: The elements
of $\t^1$ are exactly the equivalence classes
in $\tree$ of the reals in $[\alpha,\beta]$, and in particular
$\rho^1\simeq \alpha\simeq\beta$.
In a similar way, using a time reversal argument, we can construct
a subtree $\t^2$ from the right side of $\llbracket a,b\rrbracket$, with
root $\rho^2\in\rrbracket a,b \llbracket$, such that
$d_{\ov{\bf e}}(a,\rho^2)<\eta$  and
$$\inf_{c\in\t^2} \ov Z_c <\ov Z_a.$$

We can always choose $\t^1$ and $\t^2$ in such a way
that $\rho^1\prec \rho^2$. From our choice of $\eta$, we 
have then
$$\inf_{c\in\llbracket \rho^1,\rho^2\rrbracket}\ov Z_c > \ov Z_a.$$

We now exploit the convergence (\ref{basic}) to get similar 
properties for the discrete trees $(\tau^\circ_n,\ell^n)$. We can
find a positive number $\kappa$ such that the following holds for
$n$ sufficiently large. There exists a subtree $\tau^1_n$
from the left side of $\llbracket a_n,b_n\rrbracket$ with
root $\rho^1_n\in \rrbracket a_n,b_n\llbracket$ such that
\be
\label{altech1}
\inf_{x\in\tau^1_n} \ell^n_x \leq \ell^n_{a_n} -\kappa n^{1/4}.
\ee
The subtree $\tau^1_n$ is coded by an interval $[\alpha_n,\beta_n]\cap\Z$
(via the identification
$\tau^\circ_n=[pn]/\sim_n$) such that $\alpha_n/(pn)\la \alpha$
and $\beta_n/(pn)\la \beta$.
Similarly, there exists a subtree $\tau^2_n$
from the right side of $\llbracket a_n,b_n\rrbracket$ with
root $\rho^2_n\in \rrbracket a_n,b_n\llbracket$ such that
\be
\label{altech2}
\inf_{x\in\tau^2_n} \ell^n_x \leq \ell^n_{a_n} -\kappa n^{1/4}.
\ee
Furthermore, $\rho^1_n\prec \rho^2_n$ and
\be
\label{altech3}
\inf_{x\in\llbracket \rho^1_n,\rho^2_n\rrbracket} \ell^n_x \geq \ell^n_{a_n}
+\kappa n^{1/4}.
\ee

Let $\gamma_n=(\gamma_n(i),0\leq i\leq d_n(a_n,b_n))$ be a geodesic path
from $a_n$ to $b_n$ in $M_n$. As in the proof of
Lemma \ref{linesegment}, we know that the path $\gamma_n$
lies in $\tau^\circ_n$ when $n$ is large. Denote by $u_n$ the last point on the geodesic
$\gamma_n$
that does not belong to the set 
$$\{x\in\tau^\circ_n: \rho^1_n\leq x\ll \rho^2_n\}.$$
This definition makes sense because $a_n<\rho^1_n$. Also $u_n\ne b_n$ since
$\rho^1_n\leq b_n\ll \rho^2_n$. Denote by $v_n$ the point following $u_n$ on the
geodesic 
$\gamma_n$. 

Since $n^{-1/4}d_n(a_n,b_n)\la 0$ as $n\to\infty$, we know that
$$n^{-1/4}\sup_{0\leq i\leq d_n(a_n,b_n)} d_n(a_n,\gamma_n(i))
\build{\la}_{n\to\infty}^{} 0,$$
and therefore
\be
\label{altech4}
n^{-1/4}\sup_{0\leq i\leq d_n(a_n,b_n)} |\ell^n_{a_n} - \ell^n_{\gamma_n(i)} |
\build{\la}_{n\to\infty}^{} 0.
\ee

The preceding properties imply that $v_n\in \tau^1_n\cup\tau^2_n$
for $n$ sufficiently large.  Indeed, we have 
$\rho^1_n\leq v_n\ll \rho^2_n$ by construction and we also know that
$\ell^n_{v_n}>\ell^n_{a_n}-\kappa n^{1/4}$ if $n$ is large, by (\ref{altech4}).
Suppose that
$v_n\notin \tau^1_n\cup\tau^2_n$. Then, by
(\ref{altech1}) and (\ref{altech2}),
$v_n$ can
be connected to a point $y$ that does not belong to $\{x\in\tau^\circ_n:
\rho^1_n\leq x\ll
\rho^2_n\}$ only if $y\in
\llbracket \rho^1_n,\rho^2_n\rrbracket$.
Thus we get that $u_n\in
\llbracket \rho^1_n,\rho^2_n\rrbracket$, but this is impossible by
(\ref{altech3}) and (\ref{altech4}), if $n$ is large enough.

So, for $n$ sufficiently large, we have either
$v_n\in \tau^1_n$ or
$v_n\in \tau^2_n$.
One of these two cases has to occur infinitely often.
For definiteness, we assume that 
the property $v_n\in \tau^1_n$ occurs infinitely often 
and from now on until the final part of the proof we restrict our attention to integers $n$ such that
this property holds. 

Then the following properties hold for 
$n$ large:
\begin{description}
\item{\rm(i)}
$v_n\in \tau^1_n$ and $v_n\ne \rho^1_n$.
\item{\rm(ii)} $u_n\leq \rho^1_n$.
\item{\rm (iii)} Every point $w$ that comes after $v_n$ on the
geodesic $\gamma_n$ satisfies $\rho^1_n\leq w$.
\end{description}

The property $v_n\ne\rho^1_n$ is clear from (\ref{altech4}) and
(\ref{altech3}). To get (ii), recall that by construction 
we have either $u_n\leq \rho^1_n$ or $\rho^2_n\ll u_n$ (or both together).
Suppose that $\rho^2_n\ll u_n$.
If $n$ is large, the fact
that $u_n$ is connected with a point of 
$\tau^1_n\backslash\{\rho^1_n\}$ and the property (\ref{altech2}) 
then imply that  $u_n\in\llbracket \varnothing,\rho^2_n\rrbracket$.
However $u_n\in \llbracket \rho^1_n,\rho^2_n\rrbracket$ is excluded
by (\ref{altech3}) and (\ref{altech4}), and thus we get 
$u_n\in\llbracket \varnothing,\rho^1_n\rrbracket$, so that 
in particular $u_n\leq \rho^1_n$. Finally, (iii) is clear from the definition
of $u_n$.

Thanks to (i)--(iii), we can apply Lemma \ref{combi}, and we
get that if $n$ is large enough, for every point 
$y$ of $\tau^1_n\backslash\{\rho^1_n\}$ such that
\be
\label{claim0}
\ell^n_x>
\sup_{0\leq i\leq d_n(a_n,b_n)} \ell^n_{\gamma_n(i)}\;
\hbox{ for every }
x\in\llbracket\rho^1_n,y\rrbracket
\ee
we have
\be
\label{claim1}
d_n(y,b_n)\leq d_n(u_n,b_n)+\ell^n_y -
\inf_{0\leq i\leq d_n(a_n,b_n)} \ell^n_{\gamma_n(i)}.
\ee

For every $\varepsilon>0$, denote by $\u^\varepsilon_n$
the set of all vertices $y\in\tau^1_n$ such that:
\begin{description}
\item{$\bullet$} $\ell^n_y\leq \ell^n_{a_n}+ {3\varepsilon\over 2} n^{1/4}$;
\item{$\bullet$} $\ell^n_x\geq \ell^n_{a_n}+{\varepsilon\over 16} n^{1/4}$,
for every $x\in\llbracket \rho^1_n,y\rrbracket$.
\end{description}
Recall that $[\alpha,\beta]$ is the interval coding $\t^1$ 
and that $s\simeq a$.
We denote by $\u^\varepsilon_\infty$ the set of all
$r\in[\alpha,\beta]$ such that
\begin{description}
\item{$\bullet$} $\ov Z_r< \ov Z_s+
(\frac{9}{4p(p-1)})^{1/4}\,\varepsilon$;
\item{$\bullet$} $\ov Z_{r'}>\ov
Z_s+(\frac{9}{4p(p-1)})^{1/4}\,\frac{\varepsilon}{8}$, for every
$r'\in\llbracket \rho^1,r\rrbracket$.
\end{description}
(When writing $\llbracket \rho^1,r\rrbracket$ we slightly abuse notation by identifying $r$
with the corresponding vertex in $\t_{\ov{\bf e}}$.) Notice that
$\u^\varepsilon_\infty$ is open. 

Moreover, let
$]u,v[$ be a connected component of $\u^\varepsilon_\infty$,
and let $[u',v']$ be a compact subinterval of $]u,v[$. We claim that
for every $n$ sufficiently large, we must have
\be
\label{discreconti}
[pnu',pnv']\cap\Z\subset \u^\varepsilon_n
\ee
in the sense that every vertex $y$ of $\tau^\circ_n$
such that $y\sim_n k$ for some $k\in [pnu',pnv']\cap\Z$
belongs to $\u^\varepsilon_n$. To see this, first note that
the property $[pnu',pnv']\cap\Z\subset\tau^1_n$ holds for $n$
sufficiently large because $]u,v[\subset[\alpha,\beta]$. Then suppose
that for every $n$ belonging to a subsequence converging to $\infty$
we can find a vertex $y_n\in\tau^1_n$ such that 
$y_n\sim_n k_n$ for some $k_n\in [pnu',pnv']\cap\Z$
and at least one of the two conditions
\begin{description}
\item{(a)} $\ell^n_{y_n}\leq \ell^n_{a_n}+ {3\varepsilon\over 2} n^{1/4}$,
\item{(b)} $\ell^n_{x}\geq \ell^n_{a_n}+{\varepsilon\over 16} n^{1/4}$,
for every $x\in\llbracket \rho^1_n,y_n\rrbracket$,
\end{description}
does not hold. By compactness we can assume
that $k_n/(pn)\la r\in[u',v']$. 
If condition (a) fails for
infinitely many values of $n$, (\ref{basic}) gives
$$\ov Z_r\geq \ov Z_s+
\Big(\frac{9}{4p(p-1)}\Big)^{1/4}\frac{3\varepsilon}{2},$$
which contradicts the fact that $[u',v']\subset \u^\varepsilon_\infty$.
If (b) fails for infinitely values of $n$, then for these values
of $n$ we can find $\ov k_n\in[\alpha_n,k_n]\cap \Z$
such that
$$C^n_{\ov k_n}=\inf_{\ov k_n\leq k\leq k_n} C^n_k$$
and
$$V^n_{\ov k_n}< \ell^n_{a_n}+{\varepsilon\over 16} n^{1/4}.$$
Again by compactness, we can assume that $\ov k_n/(pn)\la \ov r\in[\alpha,r]$.
We have then
$$\ov{\bf e}_{\ov r}=\inf_{\ov r\leq r'\leq r} \ov{\bf e}_{r'}$$
so that $\ov r\in \llbracket \rho^1,r\rrbracket$, and
$$\ov Z_{\ov r}\leq \ov Z_s+\Big(\frac{9}{4p(p-1)}\Big)^{1/4}
\frac{\varepsilon}{16}$$
thus contradicting the fact that $[u',v']\subset \u^\varepsilon_\infty$.
This completes the proof of our claim (\ref{discreconti}).

If $I$ is a finite union of closed
subintervals of $[0,1]$ the number of vertices of 
$\tau^\circ_n=[pn]/\sim_n$ for which the first representative in $[pn]$ belongs 
to $pnI$ behaves like $(p-1)n|I|$ as $n\to\infty$, where
$|I|$ denotes the Lebesgue measure of $I$. When $I$ is of the type $[0,t]$, this
was observed in the proof of Proposition \ref{finite-margi}, and the general case
follows by a simple argument. Thus (\ref{discreconti}) implies that
\be
\label{altech5}
\liminf_{n\to\infty} \frac{1}{(p-1)n}\,\#\u^\varepsilon_n
\geq \lambda(\u^\varepsilon_\infty).
\ee

We can now use (\ref{distzero}), (\ref{altech4}) and (\ref{claim1})
to see that for $n$ sufficiently large, for every $y\in\u^\varepsilon_n
\backslash\{\rho^1_n\}$,
we have
$$d_n(y,b_n)\leq 2\varepsilon\,n^{1/4}.$$
(Notice that condition (\ref{claim0}) is satisfied for every
$y\in\u^\varepsilon_n$ when $n$ is large enough.) Hence, for every
$y,y'\in\u^\varepsilon_n$ we have also
$$d_n(y,y')\leq 4\varepsilon\,n^{1/4}.$$
Recall that $B_n(y,R)$ denotes the closed ball with radius $R$ centered
at $y$ in the metric space $({\bf m}_n,d_n)$. We have thus
$\# B_n(y,4\varepsilon n^{1/4})\geq \#\u^\varepsilon_n$
for every $y\in \u^\varepsilon_n$.

Let $(\varepsilon_k)$ be any fixed sequence
monotonically decreasing to $0$. By Lemma \ref{occuptree},
we can find $\delta_0>0$ and an integer $k_0$ such that for
every $k\geq k_0$,
$$\lambda(\u^{\varepsilon_k}_\infty) \geq 2\delta_0\,\varepsilon_k^2.$$
From (\ref{altech5}) we then see that for every $k\geq k_0$,
if $n$ is sufficiently large, we have
$$\#\u^{\varepsilon_k}_n\geq \delta_0\,\varepsilon_k^2\,(p-1)n.$$
By preceding remarks, this entails that for every $k\geq k_0$, 
if $n$ is sufficiently large,
\be 
\label{altech6}
\#\{y\in {\bf m}_n:\# B_n(y,4\varepsilon_kn^{1/4})\geq
\delta_0\varepsilon_k^2(p-1)n \} \geq \delta_0\varepsilon_k^2 (p-1)n.
\ee
Since we restricted our attention to integers $n$ such that $v_n\in\tau^1_n$,
the bound (\ref{altech6}) only holds for those integers. However, a symmetric
argument shows that (\ref{altech6}) also holds for all (sufficiently large) integers $n$
such that $v_n\in\tau^2_n$, possibly with different values of $\delta_0$ and $k_0$.
Thus by changing $\delta_0$ and $k_0$ if necessary, we can assume that 
(\ref{altech6}) holds for all sufficiently large integers $n$.

On the other hand, (\ref{pointed}) shows that for every $\delta>0$ and every $k$,
\ba
&&E\Big[\frac{1}{(p-1)n}\#\{y\in {\bf m}_n:\# B_n(y,4\varepsilon_kn^{1/4})\geq
\delta\varepsilon_k^2(p-1)n \} \Big]\\
&&\qquad\build{\la}_{n\to\infty}^{} P\Big[\ov{\cal I}\Big(\Big[0,
4\Big(\frac{9}{4p(p-1)}\Big)^{1/4}
\varepsilon_k\Big]\Big)\geq \delta\varepsilon_k^2\Big].
\ea
By Lemma \ref{estimateISE}, we have for every $\delta>0$,
$$P\Big[\ov{\cal I}\Big(\Big[0,
4\Big(\frac{9}{4p(p-1)}\Big)^{1/4}
\varepsilon_k\Big]\Big)\geq \delta\varepsilon_k^2\Big]
=o(\varepsilon_k^2)$$
as $k\to\infty$. Hence, Fatou's lemma gives
$$E\Big[\liminf_{n\to\infty}
\frac{1}{n}\#\{y\in {\bf m}_n:\# B_n(y,4\varepsilon_kn^{1/4})\geq
\delta\varepsilon_k^2(p-1)n \}\Big]=o(\varepsilon_k^2)$$
as $k\to\infty$.
Another application of Fatou's lemma yields that
$$E\Big[\liminf_{k\to\infty}\Big(
\liminf_{n\to\infty}
\frac{1}{\varepsilon_k^2n}\#\{y\in {\bf m}_n:\#
B_n(y,4\varepsilon_kn^{1/4})\geq
\delta\varepsilon_k^2(p-1)n \}\Big)\Big]=0.$$
By applying the above to a sequence of values of $\delta$ decreasing to $0$, we obtain that a.s. for every
$\delta>0$,
$$\liminf_{k\to\infty}\Big(
\liminf_{n\to\infty}
\frac{1}{\varepsilon_k^2n}\#\{y\in {\bf m}_n:\#
B_n(y,4\varepsilon_kn^{1/4})\geq
\delta\varepsilon_k^2(p-1)n \}\Big)=0.$$
This contradicts (\ref{altech6}), thus completing the proof of 
Proposition \ref{ancestorline}. \cq

\begin{proposition}
\label{notancestor}
Almost surely, for every pair $(a,b)$ in $\tree$ such that
$a$ is not an ancestor of $b$ and $b$ is not an ancestor
of $a$, the condition
$D(a,b)=0$ implies $D^\circ(a,b)=0$. 
\end{proposition}

\proof The proof is similar to that of Proposition 
\ref{ancestorline} but the fact that we already know the
property stated in this proposition makes the argument
a little simpler. We again argue by contradiction, 
assuming that there exists a pair $(a,b)$ satisfying the
condition of the proposition, such that $D(a,b)=0$
and $D^\circ(a,b)>0$. Without loss of generality we may
and will assume that $a<b$. Recall that we have
automatically $\ov Z_a=\ov Z_b$.

Let $[s,t]$ be the smallest subinterval of $[0,1]$
such that $a\simeq s$ and $b\simeq t$. As in the proof of Proposition
\ref{ancestorline}, we can find $a_n,b_n\in\tau^\circ_n$
in such a way that there exist $i_n,j_n\in [pn]$ with
$i_n\leq j_n$, $a_n\sim_n i_n$, $b_n\sim_n j_n$ and
$i_n/(pn)\la s$, $j_n/(pn)\la t$ as $n\to\infty$. We have then
$$n^{-1/4}d_n(a_n,b_n)\build{\la}_{n\to\infty}^{} D(a,b)=0.$$

From Lemma \ref{linesegment}, we also know that
$\ov Z_c\geq \ov Z_a$ for every $c\in\llbracket a,b \rrbracket$.
Recall that we assumed
$$\ov Z_a + \ov Z_b -2 \inf_{c\in[a,b]} \ov Z_c =D^\circ(a,b) >0.$$
It follows that
$$\inf_{c\in[a,b]\backslash\llbracket a,b \rrbracket} \ov Z_c
<\ov Z_a=\ov Z_b.$$
Since the minimum of $\ov {\bf e}$ over $[s,t]$ is attained at a unique
time corresponding to the vertex $a\wedge b$ (otherwise the tree $\t_{\ov{\bf e}}$ would have a point
with multiplicity strictly greater than $3$), we have
$$[a,b]\backslash\llbracket a,b \rrbracket
=([a,a\wedge b]\backslash\llbracket a,a\wedge b \rrbracket)
\cup ([a\wedge b,b]\backslash\llbracket a\wedge b,b\rrbracket).$$
Thus at least one of the following two conditions holds:
\be
\label{nota1}
\inf_{c\in[a,a\wedge b]\backslash\llbracket a,a\wedge b \rrbracket} \ov
Z_c <\ov Z_a,
\ee
or
\be
\label{nota2}
\inf_{c\in[a\wedge b,b]\backslash\llbracket a\wedge b,b\rrbracket} \ov
Z_c <\ov Z_a.
\ee
For definiteness, we assume that (\ref{nota2}) holds. The symmetric case
where (\ref{nota1}) holds is treated in a similar manner.

Under (\ref{nota2}), there exists a subtree $\t^1$ from the left
side of $\llbracket a\wedge b,b\rrbracket$, with root
$\rho^1\in \rrbracket a\wedge b,b\llbracket$, such that
$$\inf_{c\in\t^1} \ov Z_c <\ov Z_a.$$
We let $[\alpha,\beta]$ be the interval coding $\t^1$. 

As in the proof
of Proposition \ref{ancestorline}, we can find a positive number
$\kappa$ such that the following holds for $n$ sufficiently large. There exists a 
subtree $\tau^1_n$ of $\tau^\circ_n$, from the left side of
$\llbracket a_n\wedge b_n,b_n\rrbracket$, with root
$\rho^1_n\in \rrbracket a_n\wedge b_n,b_n\llbracket$
and such that
\be
\label{nota3}
\inf_{x\in \tau^1_n} \ell^n_{x} \leq \ell^n_{a_n} - \kappa n^{1/4}.
\ee
Furthermore, $\tau^1_n$ is coded by an interval $[\alpha_n,\beta_n]
\cap \Z$, with $i_n\leq \alpha_n\leq \beta_n\leq j_n$, and 
$\alpha_n/(pn)\la \alpha$, $\beta_n/(pn)\la \beta$
as $n\to\infty$. 

Let $\gamma_n=(\gamma_n(i),0\leq i\leq d_n(a_n,b_n))$ be a geodesic path
from $a_n$ to $b_n$. As previously, we know that $\gamma_n$ lies entirely in
$\tau^\circ_n$ when $n$ is large. Furthermore, as in the proof of Proposition
\ref{ancestorline}, we have
\be
\label{notatech1}
n^{-1/4}\sup_{0\leq i\leq d_n(a_n,b_n)} | \ell^n_{a_n} - \ell^n_{\gamma_n(i)} |
\build{\la}_{n\to\infty}^{} 0.
\ee
 We first
claim that for
$n$ sufficiently large the path $\gamma_n$ does not intersect $\llbracket
\varnothing,\rho^1_n
\rrbracket$. Indeed, suppose that $\gamma_n$ intersects
$\llbracket \varnothing,\rho^1_n
\rrbracket$ for infinitely many values of $n$, and for such values
write $g_n$ for one of the intersection points. Let 
$k_n\in[pn]$ be such that $g_n\sim_n k_n$, and let 
$r$ be any accumulation point of $k_n/(pn)$ in $[0,1]$. If
$c\in\tree$ is such that $c\simeq r$, the property
$d_n(g_n,b_n)\leq d_n(a_n,b_n)=o(n^{1/4})$ ensures that
$D(c,b)=0$. However, the fact that $g_n\in \llbracket \varnothing,\rho^1_n
\rrbracket$ easily implies that $c\in\llbracket \ov\rho,\rho^1\rrbracket$.
Hence we have both $D(c,b)=0$ and $c\in\llbracket \ov\rho,b\rrbracket$
with $c\ne b$. By Proposition \ref{ancestorline} this cannot occur.

Now let $u_n$ be the last point on the geodesic $\gamma_n$
that belongs to $\{x\in\tau^\circ_n:x<\rho^1_n\}$. This makes
sense since $a_n$  belongs to the latter set. Also $u_n\ne b_n$
since $\rho^1_n\leq b_n$. Let $v_n$ be the point following $u_n$
on the geodesic $\gamma_n$. We claim that $v_n\in\tau^1_n$
if $n$ is sufficiently large. Indeed, the property (\ref{nota3})
warrants that a vertex $y$ belonging to the set
$$\{x\in\tau^\circ_n: \rho^1_n\leq x\}\backslash \tau^1_n$$
and such that $\ell^n_y > \ell^n_{a_n} -\kappa n^{1/4}$ cannot be connected
to $u_n$, except possibly if
$u_n\in\llbracket \varnothing,\rho^1_n\rrbracket$. However we just saw that
this case does not occur for $n$ sufficiently large. By applying 
the preceding considerations to $y=v_n$, using (\ref{notatech1}), we get our claim.

Then the following properties hold for 
$n$ large:
\begin{description}
\item{\rm(i)}
$v_n\in \tau^1_n$ and $v_n\ne \rho^1_n$.
\item{\rm(ii)} $u_n\leq \rho^1_n$.
\item{\rm (iii)} Every point $w$ that comes after $v_n$ on the
geodesic $\gamma_n$ satisfies $\rho^1_n\leq w$.
\end{description}
From Lemma \ref{combi},
we
get that if $n$ is large enough, for every point 
$y$ of $\tau^1_n\backslash\{\rho^1_n\}$ such that
$$\ell^n_x>
\sup_{0\leq i\leq d_n(a_n,b_n)} \ell^n_{\gamma_n(i)}\;
\hbox{ for every }
x\in\llbracket\rho^1_n,y\rrbracket
$$
we have
$$d_n(y,b_n)\leq d_n(u_n,b_n)+\ell^n_y -
\inf_{0\leq i\leq d_n(a_n,b_n)} \ell^n_{\gamma_n(i)}.
$$
The end of the argument is now entirely similar to the end
of the proof of Proposition \ref{ancestorline}: We use 
(\ref{basic}), Lemma \ref{estimateISE} and Lemma \ref{occuptree} to show
that the preceding properties lead to a contradiction. This
completes the proof of Proposition \ref{notancestor}. \cq

The implication (i)$\Rightarrow$(iii) in Theorem \ref{main}
is a consequence of Propositions \ref{ancestorline}
and \ref{notancestor}. This completes the proof
of Theorem \ref{main}.

\section{Proof of the technical estimates}

In this section, we prove the three lemmas that were stated
at the end of subsection 2.4. We first need to recall some
basic properties of the Brownian snake. More information
can be found in the monograph \cite{Zu}.

The (one-dimensional) Brownian snake is a Markov process taking values
in the space $\W$ of finite paths in $\R$. Here a finite path is simply 
a continuous mapping $\w:[0,\zeta]\la \R$, where
$\zeta=\zeta_{(\w)}$ is a nonnegative real number called the 
lifetime of $\w$. The set $\W$ is a Polish space when equipped with the
distance
$$d(\w,\w')=|\zeta_{(\w)}-\zeta_{(\w')}|+\sup_{t\geq 0}|\w(t\wedge
\zeta_{(\w)})-\w'(t\wedge\zeta_{(\w')})|.$$
The endpoint (or tip) of the path $\w$ is denoted by $\wh \w=\w(\zeta_{(\w)})$.

Let
$\Omega:=C(\R_+,\W)$ be the space of all continuous functions from
$\R_+$ into $\W$, which is equipped with the topology
of uniform convergence on every compact subset of $\R_+$. 
The canonical process on $\Omega$ is then denoted by
$W_s(\omega)=\omega(s)$ for $\omega\in\Omega\;,$
and we write $\zeta_s=\zeta_{(W_s)}$ for the lifetime of $W_s$.

Let $\w\in\W$. The law of the Brownian snake started from $\w$
is the probability measure $\P_\w$ on $\Omega$ which can be 
characterized as follows. First, the process $(\zeta_s)_{s\geq 0}$
is under $\P_\w$ a reflected  Brownian motion in $[0,\infty[$ started 
from $\zeta_{(\w)}$. Secondly, the conditional distribution 
of $(W_s)_{s\geq 0}$ knowing $(\zeta_{s})_{s\geq 0}$, which is
denoted by $\Theta^\zeta_\w$, is characterized by the
following properties:
\begin{description}
\item{(i)} $W_0=\w$, $\Theta^\zeta_\w$ a.s.
\item{(ii)} The process $(W_s)_{s\geq 0}$ is time-inhomogeneous
Markov under $\Theta^\zeta_\w$. Moreover, if $0\leq s\leq
s'$,
\begin{description}
\item{$\bullet$} $W_{s'}(t)=W_{s}(t)$ for every $t\leq m(s,s'):=
\inf_{[s,s']}\zeta_r$, \ $\Theta^\zeta_\w$ a.s.
\item{$\bullet$} Under
$\Theta^\zeta_\w$, $(W_{s'}(m(s,s')+t)-W_{s'}(m(s,s')))_{0\leq t\leq
\zeta_{s'}- m(s,s')}$ is independent of $W_s$ and distributed as a
one-dimensional Brownian motion started at $0$.
\end{description}
\end{description}
Informally, the value $W_s$ of the Brownian snake at time $s$
is a random path with a random lifetime $\zeta_s$ evolving like
reflecting Brownian motion in $[0,\infty[$. When $\zeta_s$ decreases,
the path is erased from its tip, and when $\zeta_s$ increases, the path 
is extended by adding ``little pieces'' of Brownian paths at its tip.

We denote by $n(de)$ the It\^o measure of positive Brownian excursions,
which is a $\sigma$-finite measure on the space $C(\R_+,\R_+)$,
and we write
$$\sigma(e)=\inf\{s> 0:e(s)=0\}$$
for the duration of excursion $e$. 
For $s>0$, $n_{(s)}$ will denote the conditioned measure $n(\cdot\mid\sigma=s)$. In particular
$n_{(1)}(de)$ is the law of the normalized excursion $\eg$, or more 
precisely of $(\eg_{t\wedge 1})_{t\geq 0}$.
Our normalization of the excursion measure is fixed by the relation
\be
\label{decomItomeas}
n=\int_0^\infty {ds\over 2\sqrt{2\pi s^3}}\;n_{(s)},
\ee
and we have then $n(\sup_{s\geq 0}e(s)>\varepsilon)=(2\varepsilon)^{-1}$ for every $\varepsilon>0$.

 If $x\in\R$, the excursion measure
$\N_x$ of the Brownian snake from $x$ is given by
$$\N_x=\int_{C(\R_+,\R_+)} n(de)\;\Theta^e_{\ov x}$$
where $\ov x$ denotes the trivial element of $\W$ with lifetime $0$
and initial point $x$. With a slight abuse of notation
we also write $\sigma(\omega)=\inf\{s>0:\zeta_s(\omega)=0\}$ for
$\omega\in\Omega$. We can then consider the conditioned measures
$$\N_x^{(s)}=\N_x(\cdot\mid \sigma=s)=
\int_{C(\R_+,\R_+)} n_{(s)}(de)\;\Theta^e_{\ov x}\;.$$

We can now relate the Brownian snake to the Brownian trees of
subsection 2.4: We may define the pair $(\eg,Z)$ under the probability measure 
$\N^{(1)}_0$ by taking $\eg_s=\zeta_s$ and $Z_s=\wh W_s$, for every $0\leq s\leq 1$.
Furthermore the path $(W_s(t),0\leq t\leq \zeta_s)$ is then interpreted in terms
of the labels attached to the ancestors of $p_\eg(s)$: If $a=p_\eg(s)$ is a vertex of the
tree
$\t_\eg$, and
$c\in\llbracket \rho,a\rrbracket$ is the ancestor of $a$ at generation $t=d_\eg(\rho,c)$,
we have $Z_c=W_s(t)$. These identifications follow very easily from the
properties of the Brownian snake.

For future reference, we state a crude 
bound on the increments of the process $(\wh W_s)_{s\geq 0}$ under $\N_0$. In a way analogous
to subsection 2.3 we set, for every $s,t\geq 0$,
$$d_\zeta(s,t)=\zeta_s+\zeta_t-2\inf_{s\wedge t\leq r\leq s\vee t} \zeta_r\;.$$

\begin{lemma}
\label{crudebound}
Let $b\in]0,1/2[$. Then $\N_0(d\omega)$ a.e. there exists $\varepsilon_0(\omega)>0$ such that
for every $s,t\geq 0$ with $d_\zeta(s,t)\leq \varepsilon_0$, one has
$$|\wh W_s-\wh W_t|\leq (d_\zeta(s,t))^b.$$
\end{lemma}

\proof Conditionally on $(\zeta_r)_{r\geq 0}$, the process $(\wh W_r)_{r\geq 0}$ is
Gaussian with mean $0$ and such that $E[(\wh W_s-\wh W_t)^2]=d_\zeta(s,t)$
for every $s,t\geq 0$. The bound of the lemma then follows
from standard chaining arguments. We leave details to the reader. \cq

\smallskip
\noindent{\bf Proof of Lemma \ref{estimateISE}:} Set
$$\un W=\inf_{s\geq 0} \wh W_s,$$
and, for every $\varepsilon>0$, 
$${\cal J}(\varepsilon)=\int_0^\sigma ds\,{\bf 1}{\{\wh W_s-\un W\leq \varepsilon\}}.$$
Thanks to the remarks preceding Lemma \ref{crudebound}, the quantity $\ov{\cal I}([0,\varepsilon])$
in Lemma \ref{estimateISE} has the same distribution as ${\cal J}(\varepsilon)$
under $\N^{(1)}_0$. Therefore, the
statement of Lemma \ref{estimateISE} reduces to checking that
$$\N^{(1)}_0({\cal J}(\varepsilon)\geq \alpha \varepsilon^2)=o(\varepsilon^2)$$
as $\varepsilon\to 0$. From (\ref{decomItomeas}) and simple scaling arguments, it is enough to verify
that
\be
\label{estmISEtech}
\N_0({\cal J}(\varepsilon)\geq \alpha \varepsilon^2\,,\,\sigma>1/2)=o(\varepsilon^2)
\ee
as $\varepsilon\to 0$. For every $\delta>0$, we have
${\cal J}(\varepsilon)={\cal J}_\delta(\varepsilon)+{\cal J}'_\delta(\varepsilon)$,
where 
$${\cal J}_\delta(\varepsilon)=\int_0^\sigma ds\,{\bf 1}{\{\wh W_s-\un W\leq \varepsilon,
\zeta_s<\delta\}}\ ,\ 
{\cal J}'_\delta(\varepsilon)=\int_0^\sigma ds\,{\bf 1}{\{\wh W_s-\un W\leq \varepsilon,
\zeta_s\geq\delta\}}.
$$
Let us fix
$\beta>0$. By Lemma 3.2 in \cite{LGW}, we can choose $\delta>0$ small enough so that,
for every $\varepsilon\in]0,1[$,
$$\N_0\Big({\bf 1}_{\{\sigma >1/2\}}\,{\cal J}_\delta(\varepsilon)\Big)\leq \beta\,\varepsilon^4.$$
On the other hand, Lemma 3.3 in \cite{LGW} yields the existence of a constant $K_\delta$
such that, for every $\varepsilon\in]0,1[$,
$$\N_0\Big(({\cal J}'_\delta(\varepsilon))^2\Big)\leq K_\delta\,\varepsilon^8.$$
Then,
\ba
\N_0({\cal J}(\varepsilon)\geq \alpha \varepsilon^2\,,\,\sigma>1/2)&\leq&
\N_0({\cal J}_\delta(\varepsilon)\geq \frac{\alpha}{2}\,\varepsilon^2\,,\,\sigma>1/2)
+\N_0({\cal J}'_\delta(\varepsilon)\geq \frac{\alpha}{2}\,\varepsilon^2)\\
&\leq& \frac{2}{\alpha \varepsilon^2}
\N_0\Big({\bf 1}_{\{\sigma >1/2\}}\,{\cal J}_\delta(\varepsilon)\Big)
+\frac{4}{\alpha^2\varepsilon^4}\,\N_0\Big(({\cal J}'_\delta(\varepsilon))^2\Big)\\
&\leq& \frac{2\beta}{\alpha}\,\varepsilon^2+\frac{4 K_\delta}{\alpha^2}\,\varepsilon^4.
\ea
It follows that
$$\limsup_{\varepsilon\to 0} \varepsilon^{-2}\N_0({\cal J}(\varepsilon)\geq \alpha
\varepsilon^2\,,\,\sigma>1/2)\leq \frac{2\beta}{\alpha}.$$
Since $\beta$ was arbitrary, this completes the proof of (\ref{estmISEtech}) and of the lemma. \cq

\medskip
\noindent{\bf Proof of Lemma \ref{increasepoint}:}
We first explain why it is enough to prove the statement concerning
the pair $(\eg,Z)$. This follows from a re-rooting argument. Recall
the notation of subsection 2.4. For every fixed $s\in[0,1[$, set
\begin{description}
\item{$\bullet$} 
$\displaystyle{\eg^{[s]}_t=\eg_{s}+\eg_{s\oplus t}-2\,m_{\eg}(s,
{s\oplus t})}$;
\item{$\bullet$} $Z^{[s]}_t=Z_{s\oplus t} -Z_{s}$,
\end{description}
for every $t\in[0,1]$. By construction, $(\ov\eg,\ov Z)=(\eg^{[s_*]},Z^{[s_*]})$. 
Also, $(\eg^{[s]},Z^{[s]})\build{=}_{}^{\rm(d)} (\eg,Z)$ for every fixed
$s\in[0,1[$: See Proposition 4.9 in \cite{MaMo} or
Theorem 2.3 in \cite{LGW}. Hence, if $U$ is uniformly distributed
over $[0,1[$ and independent of $(\eg,Z)$, we have also 
$(\eg^{[U]},Z^{[U]})\build{=}_{}^{\rm(d)} (\eg,Z)$.

Suppose there exists an increase point $r\in]0,1[$ of the pair
$(\ov\eg,\ov Z)=(\eg^{[s_*]},Z^{[s_*]})$. Then for every
$s$ sufficiently close to $s_*$, $r+s_*-s$ will be an increase point
of the pair $(\eg^{[s]},Z^{[s]})$ (this can be verified by
direct inspection of the formulas defining the pair
$(\eg^{[s]},Z^{[s]})$, keeping in mind that $s_*$ corresponds to
a leaf of the tree $\t_{\eg}$, so that immediately after 
or immediately before $s_*$, $\eg_s$ takes values strictly less
than $\eg_{s_*}$). In particular, the pair $(\eg^{[U]},Z^{[U]})$
will have an increase point with positive probability, which contradicts the
first assertion of the lemma. 

Let us now prove the statement concerning
the pair $(\eg,Z)$. In terms of the Brownian snake, we need
to check that $\N^{(1)}_0$ a.s. the pair $(\zeta_s,\wh W_s)_{0\leq s\leq 1}$
has no increase point. By a simple scaling argument, it is enough to
verify that the same property holds for the pair $(\zeta_s,\wh W_s)_{0\leq s\leq \sigma}$
under the excursion measure $\N_0$ (obviously time $1$ is now replaced by $\sigma$
in the definition of an increase point). To this end, we will use the
following lemma.

\begin{lemma}
\label{increaselemma}
Let $\delta>0$. Let $\w\in{\cal W}$ with $\w(0)=0$ and $\zeta_{(\w)}=a>0$, and
let $\varepsilon\in]0,a]$. Consider the
stopping times 
\ba
&&T=\inf\{s\geq 0:\zeta_s=a+\delta\},\\
&&T'=\inf\{s\geq 0:\zeta_s=a-\varepsilon\}.
\ea
On the event $\{T<T'\}$, we also define
$$L=\sup\{s\leq T:\zeta_s=a\}.$$
Then there exists a constant $C_\delta$, which only depends on $\delta$, such that,
for every $\eta\in]0,1]$,
$$\P_\w(T<T'\hbox{ and }\wh W_s>\wh W_L-\eta\hbox{ for every }s\in[L,T])\leq
C_\delta\,\varepsilon\,\eta^3.$$
\end{lemma}

\noindent{\bf Remark.} The exponent $3$ in $\eta^3$ is sharp and related to the
fact that the bound of the lemma is a ``one-sided'' estimate. This should be
compared with the exponent $4$ that appears in similar two-sided estimates 
derived in \cite{LGW}.

\proof Under $\P_\w$, $(\zeta_s)_{s\geq 0}$ is distributed as a reflected linear 
Brownian motion started from $a$. In particular,
$$\P_\w(T<T')=\frac{\varepsilon}{\varepsilon+\delta}.$$
Moreover, from standard connections between linear Brownian motion and
the three-dimen\-sional Bessel process, we know that under the conditional
probability $\P_\w(\cdot\mid T<T')$, the shifted process
$$Y_s:=\zeta_{(L+s)\wedge T}-a\ ,\ s\geq 0$$
is distributed as a three-dimensional Bessel process started from $0$ and stopped
when it first hits $\delta$. At this point, it is convenient to introduce the
future infimum process of $Y$,
$$J_s:=\inf_{r\geq s} Y_r\ ,\ s\geq 0$$
and the excursions of $Y-J$ away from $0$: Let $]\alpha_i,\beta_i[$, $i\in I$, be the connected 
components of the open set $\{s\geq 0:Y_s>J_s\}$, and for every $i\in I$ set
$$e_i(s)=Y_{(\alpha_i+s)\wedge \beta_i}-Y_{\alpha_i}.$$
Then the point measure
$$\sum_{i\in I} \delta_{(Y_{\alpha_i},e_i)}(dt\,de)$$
is Poisson with intensity
$$2 \,{\bf 1}_{[0,\delta]}(t)\,dt\;n\Big(de \cap\Big\{\sup_{s\geq 0}e(s)<\delta-t\Big\}\Big).$$
The last property follows from standard facts of excursion theory. See e.g. Lemma 1
in \cite{AW} for a detailed derivation. 

We can then combine the preceding excursion decomposition of the paths of $Y$
with the spatial displacements of the Brownian snake, in a way similar
to the proof of Lemma V.5 in \cite{Zu}. Let $H=\sup_{s\geq 0}\zeta_s$
denote the maximum of the lifetime process. It follows that
\begin{eqnarray}
\label{increasetech1}
&&\P_\w(T<T'\hbox{ and }\wh W_s>\wh W_L-\eta\hbox{ for every }s\in[L,T])\nonumber\\
&&\quad = \frac{\varepsilon}{\varepsilon+\delta}\,
E_0\Big[{\bf 1}_{\{\xi[0,\delta]\subset]-\eta,\infty[\}}\,\exp\Big(-2\int_0^\delta dt
\,\N_{\xi_t}(H<\delta-t,\un W\leq -\eta)\Big)\Big]
\end{eqnarray}
where $(\xi_t)_{t\geq 0}$ is a linear Brownian motion started from $x$ under the
probability measure $P_x$, and we use the notation $\xi[0,\delta]$ for the range of
$\xi$ over the time interval $[0,\delta]$. 

From this point, the argument is very similar to the end of the proof of Proposition 4.2 in \cite{LGW},
to which we refer the reader for more details. For every
$x>0$, we set
$$f(x)=\N_0(\un W>-x\mid H=1)$$
and
$$G(x)=4\int_0^x u(1-f(u))\,du.$$
Note that $G(+\infty)=6$ (see Section 4 in \cite{LGW}).
By conditioning with respect to $H$ and then using a scaling argument, we get
$$\int_0^\delta dt\,\N_{\xi_t}(H<\delta-t,\un W\leq -\eta)
=\int_0^\delta dt\int_0^{\delta-t} {du\over 2 u^2}\,(1-f(\frac{\xi_t+\eta}{\sqrt{u}})).$$
Hence, the right-hand side of (\ref{increasetech1}) can be written as
\ba
&&\frac{\varepsilon}{\varepsilon+\delta}\,
E_0\Big[{\bf 1}_{\{\xi[0,\delta]\subset]-\eta,\infty[\}}\,\exp\Big(-\int_0^\delta dt
\int_0^{\delta-t} {du\over u^2}\,(1-f(\frac{\xi_t+\eta}{\sqrt{u}}))\Big)\Big]\\
&&=\frac{\varepsilon}{\varepsilon+\delta}\,
E_\eta\Big[{\bf 1}_{\{\xi[0,\delta]\subset]0,\infty[\}}\,\exp\Big(-\int_0^\delta dt
\int_0^{\delta-t} {du\over u^2}\,(1-f(\frac{\xi_t}{\sqrt{u}}))\Big)\Big].
\ea
From the definition of $G$, the property $G(+\infty)=6$ and a change of variables, we have
$$\int_0^{\delta-t} {du\over u^2}\,(1-f(\frac{\xi_t}{\sqrt{u}}))
=(\xi_t)^{-2} \Big(3-\frac{1}{2}G(\frac{\xi_t}{\sqrt{\delta-t}})\Big).$$
Hence we get
\begin{eqnarray}
\label{increasetech2}
&&\P_\w(T<T'\hbox{ and }\wh W_s>\wh W_L-\eta\hbox{ for every }s\in[L,T])\nonumber\\
&&\quad =\frac{\varepsilon}{\varepsilon+\delta}\,
E_\eta\Big[{\bf 1}_{\{\xi[0,\delta]\subset]0,\infty[\}}\,\exp\Big(-3\int_0^\delta \frac{dt}{\xi_t^2}
+\frac{1}{2}\int_0^\delta \frac{dt}{\xi_t^2}\,
G(\frac{\xi_t}{\sqrt{\delta-t}})\Big)\Big].
\end{eqnarray}
Proposition 2.6 of \cite{LGW}, which reformulates absolute continuity relations between
Bessel processes due to Yor, implies that the right-hand side of (\ref{increasetech2})
is equal to
$$\frac{\varepsilon}{\varepsilon+\delta}\,\eta^3\,E^{(7)}_\eta\Big[
(R_\delta)^{-3}\,\exp\Big(\frac{1}{2}\int_0^\delta \frac{dt}{R_t^2}\,G(\frac{R_t}{\sqrt{\delta-t}})\Big)
\Big],$$
where $(R_t)_{t\geq 0}$ is a Bessel process of dimension $7$ started from $\eta$
under the probability measure $P^{(7)}_\eta$. Finally, we can argue as in the end 
of the proof of Proposition 4.2 in \cite{LGW} to verify the existence of a
constant $C'_\delta$ such that, for every $\eta>0$,
$$E^{(7)}_\eta\Big[
(R_\delta)^{-3}\,\exp\Big(\frac{1}{2}\int_0^\delta \frac{dt}{R_t^2}\,G(\frac{R_t}{\sqrt{\delta-t}})\Big)
\Big]
\leq C'_\delta.$$
Lemma \ref{increaselemma} follows with $C_\delta=\delta^{-1}C'_\delta$. \cq

\smallskip
We come back to the proof of Lemma \ref{increasepoint}. We fix $\delta\in]0,1[$. For every
$\varepsilon\in]0,1[$, we introduce the sequence of stopping times defined
inductively by
$$T^\varepsilon_0=0\ ,\
T^\varepsilon_{i+1}=\inf\{s>T^\varepsilon_i:|\zeta_s-\zeta_{T^\varepsilon_i}|=\varepsilon\},$$
with the usual convention $\inf\varnothing=\infty$. For every index $i$ such that $T^\varepsilon_i<\infty$,
we also set
\ba
&&S^\varepsilon_i=\inf\{s>T^\varepsilon_i:\zeta_s=\zeta_{T^\varepsilon_i}+\delta\},\\
&&\wt T^\varepsilon_i=\inf\{s>T^\varepsilon_i:\zeta_s=\zeta_{T^\varepsilon_i}-\varepsilon\}.
\ea
On the event $\{T^\varepsilon_i=\infty\}$ simply set $S^\varepsilon_i=\wt T^\varepsilon_i=\infty$.
Finally, on the event $\{T^\varepsilon_i<\infty\}\cap \{S^\varepsilon_i<\infty\}$, we put
$$L^\varepsilon_i=\sup\{s<S^\varepsilon_i:\zeta_s=\zeta_{T^\varepsilon_i}\}.$$
Fix $A>0$, and let ${\cal A}_{\varepsilon,\eta}$ be the event that there exists
$i\geq 1$ such that $T^\varepsilon_i<\infty$, $\zeta_{T^\varepsilon_i}\in]0,A]$, $S^\varepsilon_i<\wt
T^\varepsilon_i<\infty$ and
$$\wh W_s> \wh W_{L^\varepsilon_i}-\eta$$
for every $s\in[L^\varepsilon_i,S^\varepsilon_i]$. 

From Lemma \ref{increaselemma} and the strong Markov property for
the Brownian snake, we have
\ba
\N_0({\cal A}_{\varepsilon,\eta})
\!\!\!&\leq&\!\!\!\N_0\Big(\!\sum_{i=1}^\infty {\bf
1}{\{T^\varepsilon_i<\infty,\zeta_{T^\varepsilon_i}\in]0,A]\}}
\,{\bf 1}{\{S^\varepsilon_i<\wt
T^\varepsilon_i<\infty\}}\,{\bf 1}{\{\wh W_s>\wh W_{L^\varepsilon_i}-\eta\,,\,\forall s
\in[L^\varepsilon_i,S^\varepsilon_i]\}}\!\Big)\\
\!\!\!&\leq&\!\!\! C_\delta \varepsilon\,\eta^3\,\N_0\Big(
\sum_{i=1}^\infty {\bf 1}{\{T^\varepsilon_i<\infty,\zeta_{T^\varepsilon_i}\in]0,A]\}}\Big).
\ea
Standard properties of linear Brownian motion give
$$\N_0\Big(\sum_{i=1}^\infty {\bf 1}{\{T^\varepsilon_i<\infty,\zeta_{T^\varepsilon_i}\in]0,A]\}}\Big)
=\frac{1}{\varepsilon}\,\lfloor\frac{A}{\varepsilon}\rfloor.$$
Therefore we have obtained the bound
$$\N_0({\cal A}_{\varepsilon,\eta})\leq C_\delta A\,\varepsilon^{-1}\eta^3.$$
We apply this estimate with $\varepsilon=\varepsilon_p=2^{-p}$, for every integer $p\geq 1$,
and $\eta=(\varepsilon_p)^b$, where $b\in]\frac{1}{3},\frac{1}{2}[$. It follows that
$\N_0$ a.e. for all $p$ sufficiently large the event ${\cal A}_{\varepsilon_p,(\varepsilon_p)^b}$
does not occur. 

To complete the argument, notice that it is enough to prove that there cannot exist $r>0$ such that
$\inf\{u\geq r:\zeta_u=\zeta_r+2\delta\}<\infty$ and 
$$\zeta_s\geq \zeta_r\hbox{ and }\wh W_s\geq \wh W_r\;,\ \hbox{for every } s\in[r,\inf\{u\geq
r:\zeta_u=\zeta_r+2\delta\}].$$
We argue by contradiction and suppose that there is such a value of $r$. Let $i\geq 1$
be such that $r\in]T^{\varepsilon_p}_{i-1},T^{\varepsilon_p}_i]$. If $p$ has been taken
large enough, we have $S^{\varepsilon_p}_i<\wt T^{\varepsilon_p}_i\wedge 
\inf\{u\geq r:\zeta_u=\zeta_r+2\delta\}$, and 
for every $s\in[T^{\varepsilon_p}_i,S^{\varepsilon_p}_i]$,
$$\wh W_s\geq \wh W_r > W_{L^{\varepsilon_p}_i}-(\varepsilon_p)^b\;,$$
where the last inequality follows from Lemma \ref{crudebound} since 
$d_\zeta(r,L^{\varepsilon_p}_i)< \varepsilon_p$. We thus get a contradiction with the
fact that ${\cal A}_{\varepsilon_p,(\varepsilon_p)^b}$ does not occur when $p$
is large. This contradiction completes the proof. \cq

\medskip
\noindent{\bf Proof of Lemma \ref{occuptree}:} We first observe that it
is enough to prove  the statement of Lemma \ref{occuptree} when the pair $(\t_{\ov\eg},\ov Z)$ is replaced 
by $(\t_\eg,Z)$, and of course $\ov\rho$ is also replaced by the root $\rho$
of $\t_\eg$. This follows from a re-rooting argument analogous to the one we used at
the beginning of the proof of Lemma \ref{increasepoint}. Let us only sketch
the argument. We assume that the property of Lemma \ref{occuptree} has been derived 
when the pair $(\t_{\ov\eg},\ov Z)$ is replaced 
by $(\t_\eg,Z)$. Suppose that the
conclusion of this lemma fails for some subtree of $\t_{\ov\eg}$. Then it
will also fail for some subtree of the re-rooted tree $\t_{\eg^{[s]}}$, provided 
that $s$ is sufficiently close to $s_*$. Hence with positive probability it will
fail for some subtree of $\t_{\eg^{[U]}}$, where $U$ is uniformly distributed
over $[0,1]$. Since we saw that $(\eg^{[U]},Z^{[U]})\build{=}_{}^{\rm(d)} (\eg,Z)$,
this leads to a contradiction.

Then, we notice that by a symmetry 
argument we need only consider subtrees of $\t_\eg$ from the right side of $\llbracket\rho,a\rrbracket$.
Furthermore, as we already observed, the pair
$(\eg_s,Z_s)_{0\leq s\leq 1}$ has the same distribution as $(\zeta_s,\wh W_s)_{0\leq s\leq 1}$ under
$\N^{(1)}_0$. By scaling, it is then enough to prove that the analogue of Lemma \ref{occuptree} holds for
the pair $(\zeta_s,\wh W_s)_{0\leq s\leq \sigma}$ under $\N_0$. We can thus reformulate
the desired property in the following way. Let us fix $s>0$, and argue on the event $\{s<\sigma\}$.
Denote by $a=p_\zeta(s)$ the vertex corresponding to $s$ in the tree $\t_\zeta$.
The subtrees of $\t_\zeta$ from the right side of $\llbracket \rho,a\rrbracket$
exactly correspond to the excursions of the shifted process $(\zeta_{s+r})_{r\geq 0}$
above its past minimum process. More precisely, set
$$\zeta^{(s)}_r=\zeta_{s+r}\ ,\ \check\zeta^{(s)}_r=\inf_{0\leq u\leq r}\zeta^{(s)}_u$$
for every $r\geq 0$. Denote by $]\alpha_i,\beta_i[$, $i\in I$, the connected components
of the open set $\{r\geq 0:\zeta^{(s)}_r>\check\zeta^{(s)}_r\}$. Then for each
$i\in I$, the set $\t^1_i:=p_\zeta([s+\alpha_i,s+\beta_i])$ is a subtree of $\t_\zeta$
from the right side of $\llbracket \rho,a\rrbracket$ with root
$p_\zeta(s+\alpha_i)=p_\zeta(s+\beta_i)$, and conversely all subtrees 
from the right side of $\llbracket \rho,a\rrbracket$ are obtained in this way.
Recall the interpretation of the path $(W_{s+r}(t),0\leq t\leq \zeta_{s+r})$ as giving the
labels of the ancestors of the vertex $p_\zeta(s+r)$ in the tree $\t_\zeta$.
In order to get the statement of Lemma \ref{occuptree}, it is enough to
prove the following claim.

\smallskip
\noindent{\bf Claim}. {\it $\N_0$ a.e. on the event $\{s<\sigma\}$, for every $\mu>0$ and every
$i\in I$ such that
\be
\label{occuptech0}
\inf_{\alpha_i\leq r\leq \beta_i} \wh W_{s+r} < \wh W_{s+\alpha_i} -\mu
\ee
we have 
\be
\label{occuptech1}
\liminf_{\varepsilon\to 0}
\varepsilon^{-2}\!\! \int_{\alpha_i}^{\beta_i} dr\,
{\bf 1}{\{\wh W_{s+r}\leq \wh W_{s+\alpha_i}-\mu+\varepsilon\}}
{\bf 1}{\{W_{s+r}(t)\geq \wh W_{s+\alpha_i}-\mu+\frac{\varepsilon}{8},
\forall t\in[\zeta_{s+\alpha_i},\zeta_{s+r}]\}}>0.
\ee}

Note that the preceding claim is concerned with subtrees from the right side of one particular vertex
$a=p_\zeta(s)$, whereas the statement of Lemma \ref{occuptree} holds simultaneously for all 
choices of the vertex $a$. However, assuming that the claim is proved, it immediately follows that
the desired property holds for all subtrees from the left side of $\llbracket \rho, p_\zeta(s)\rrbracket$,
for all rational numbers $s>0$, outside a single set of zero $\N_0$-measure. Since a subtree 
from the right side of $\llbracket \rho,p_\zeta(s)\rrbracket$ is also a subtree from the right side of
$\llbracket \rho,p_\zeta(s')\rrbracket$ as soon as $s'$ is close enough to $s$, we then get the desired
result simultaneously for all choices of $a=p_\zeta(s)$. 

\smallskip
Let us now discuss the proof of the claim. Recall that $s>0$ is fixed and that we
argue on the event $\{s<\sigma\}$. For every $i\in I$ and every $r\geq 0$ set
\ba
&&\zeta^i_r=\zeta_{(s+\alpha_i+r)\wedge (s+\beta_i)}-\zeta_{s+\alpha_i}\\
&&W^i_r(t)=W_{(s+\alpha_i+r)\wedge (s+\beta_i)}(\zeta_{s+ \alpha_i}+t)-\wh W_{s+\alpha_i}\ ,\ 
\hbox{ for every }t\in[0,\zeta^i_r]
\ea
and view $W^i_r$ as a finite path with lifetime $\zeta^i_r$, so that $W^i=(W^i_r)_{r\geq 0}$
is a random element of $\Omega=C(\R_+,{\cal W})$. Also set
$\sigma_i=\beta_i-\alpha_i$, which corresponds to the duration of
the ``excursion'' $\zeta^i$. By combining the Markov property at time
$s$ with Lemma V.5 in \cite{Zu}, we get that under the probability
measure $\N_0(\cdot\mid s<\sigma)$ and conditionally on $W_s$, the point measure
$$\sum_{i\in I} \delta_{W^i}(d\omega)$$
is Poisson on $\Omega$ with intensity
$2\,\zeta_s\,\N_{0}(d\omega)$.
Now observe that condition (\ref{occuptech0}) reduces to
$$\inf_{r\geq 0} \wh W^i_r<-\mu$$
and that the integral in (\ref{occuptech1}) is equal to
$$\int_{0}^{\sigma_i} dr\,
{\bf 1}{\{\wh W^i_{r}\leq -\mu+\varepsilon\}}\,
{\bf 1}{\{W^i_{r}(t)\geq  -\mu+\frac{\varepsilon}{8}\;,\; 
\forall t\in[0,\zeta^i_{r}]\}}.$$
Thanks to these observations and to our previous description of the conditional
distribution of the point measure $\sum_{i\in I} \delta_{W^i}$, we see that our claim
follows from the next lemma.

\begin{lemma}
\label{occuplast}
$\N_0$ a.e. for every $\mu\in]0,-\un W[$, we have
$$\liminf_{\varepsilon\to 0}
\varepsilon^{-2} \int_{0}^{\sigma} dr\,
{\bf 1}{\{\wh W_{r}\leq -\mu+\varepsilon\}}\;
{\bf 1}{\{W_{r}(t)\geq  -\mu+\frac{\varepsilon}{8}\;,\; 
\forall t\in[0,\zeta_{r}]\}}>0.$$
\end{lemma}

\noindent{\bf Proof of Lemma \ref{occuplast}.} We fix an integer $N\geq 2$. Without
loss of generality, we may and will 
restrict our attention to values $\mu\in[N^{-1},N]$.
We also consider another integer $n\geq N$. If $j$ is
the integer such that $(j-1)2^{-n-3}< \mu \leq j2^{-n-3}$, and 
if $2^{-n-1}\leq \varepsilon\leq 2^{-n}$,
we have the following simple
inequalities:
\ba
&&{\bf 1}{\{\wh W_{r}\leq -\mu+\varepsilon\}}\;
{\bf 1}{\{W_{r}(t)\geq  -\mu+\frac{\varepsilon}{8}\;,\; 
\forall t\in[0,\zeta_{r}]\}}\\
&&\qquad\geq {\bf 1}{\{\wh W_{r}\leq -\mu+2^{-n-1}\}}\;
{\bf 1}{\{W_{r}(t)\geq  -\mu+2^{-n-3}\;,\; 
\forall t\in[0,\zeta_{r}]\}}\\
\noalign{\smallskip}
&&\qquad\geq {\bf 1}{\{\wh W_{r}\leq -j2^{-n-3}+2^{-n-1}\}}\;
{\bf 1}{\{W_{r}(t)\geq  -j2^{-n-3}+2^{-n-2}\;,\; 
\forall t\in[0,\zeta_{r}]\}}.
\ea

So, for every integer $j$ such that $N^{-1}2^{n+3}\leq j \leq N 2^{n+3}$, we set
$$U_{n,j} =\int_0^\sigma dr\,{\bf 1}{\{\wh W_{r}\leq -j2^{-n-3}+2^{-n-1}\}}\;
{\bf 1}{\{W_{r}(t)\geq  -j2^{-n-3}+2^{-n-2}\;,\; 
\forall t\in[0,\zeta_{r}]\}}.$$
For every $r>0$, denote by $L^r$ the total mass of the exit measure 
of the Brownian snake from
the open set $]-r,\infty[$ (see e.g. Chapter 6 of \cite{Zu} for the definition
and main properties of exit measures). Note that $\{\un W<-r\}=\{L^r>0\}$,
$\N_0$ a.e. Put
$$r_{n,j}:=-j2^{-n-3}+2^{-n-1}\leq -N^{-1}/2<0$$
to simplify notation. By the special Markov property (cf Section 2.4 in \cite{LGW}), 
conditionally on $\{L^{r_{n,j}}=\ell\}$, the variable $U_{n,j}$ is distributed as
$$\int \n(d\omega)\,X_n(\omega)$$
where $\n$ is a Poisson point measure with intensity $\ell \N_0$, and
$$X_n=\int_0^\sigma dr\,{\bf 1}{\{\wh W_{r}\leq 0\}}\;
{\bf 1}{\{W_{r}(t)\geq  -2^{-n-2}\;,\; 
\forall t\in[0,\zeta_{r}]\}}.$$
From scaling properties of $\N_0$, 
$$\int \n(d\omega)\,X_n(\omega)\build{=}_{}^{\rm(d)} 2^{-4n}\int \n_n(d\omega)\,X_0(\omega),$$
where $\n_n$ is a Poisson point measure with intensity $\ell 2^{2n}\N_0$. Note that the
quantity
$$2^{-2n}\int \n_n(d\omega)\,X_0(\omega)$$
is the mean of $2^{2n}$ independent nonnegative random variables distributed as
$\int \n(d\omega)\,X_0(\omega)$. 

We can then use standard large deviations estimates for sums of i.i.d.
random variables to derive the following. If $\eta>0$ is fixed, we can find
two positive constants $\nu$ and $\kappa$ such that, for every $n$ large enough, 
for every integer $j\geq N^{-1}2^{n+3}$, 
$$\N_0(\{2^{2n}U_{n,j}\leq \nu\}\cap\{L^{r_{n,j}}\geq \eta\})\leq \exp(-\kappa
2^{2n})\,\N_0(L^{r_{n,j}}\geq \eta)\leq c_0\,\exp(-\kappa
2^{2n}),$$
where $c_0=\N_0(\un W\leq -N^{-1}/2)$ is a positive constant. In the last
inequality we use the fact that $\N_0(L^{r_{n,j}}\geq \eta)\leq \N_0(L^{r_{n,j}}>0)
= \N_0(\un W\leq r_{n,j})$. We can sum the preceding estimate over values
of $j\in [N^{-1}2^{n+3},N 2^{n+3}]$, and then use the Borel-Cantelli lemma to
get that $\N_0$ a.e. for all $n$ sufficiently large and all $j\in [N^{-1}2^{n+3},N 2^{n+3}]$ we have either 
$L^{r_{n,j}}<\eta$ or $U_{n,j} > \nu \,2^{-2n}$.

Now recall the elementary inequalities of the beginning of the proof. It follows that
$\N_0$ a.e., for all $\mu\in[N^{-1},N]$ we have either
\begin{equation}
\label{lasttech1}
\inf_{r\in [-\mu,-(2N)^{-1}]\cap \Q} L^r \leq \eta
\ee
or, for $\varepsilon$ small enough, 
\be
\label{lasttech2}
\int_{0}^{\sigma} dr\,
{\bf 1}{\{\wh W_{r}\leq -\mu+\varepsilon\}}\;
{\bf 1}{\{W_{r}(t)\geq  -\mu+\frac{\varepsilon}{8}\;,\; 
\forall t\in[0,\zeta_{r}]\}}\geq \nu \varepsilon^2.
\ee

A simple application of the special Markov property shows that under the 
probability measure $\N_0(\cdot\mid \un W\leq (2N) ^{-1})$ the process
$(L^{-(2N) ^{-1}-a})_{a\geq 0}$ is a continuous-state branching process, hence
a Feller Markov process which is absorbed at the origin. Thus, for every $a>0$, we have
\be
\label{lasttech3}\inf_{r\in [-(2N)^{-1}-a,-(2N)^{-1}]\cap \Q} L^r >0,\qquad
\N_0\hbox{ a.e. on the event }\{\un W< -(2N)^{-1}-a\}.
\ee
We now take
$\eta=\eta_k=2^{-k}$, for every integer $k\geq 1$ (then $\nu=\nu_k$ also depends on $k$). If
$$\mu_k=\inf\{a\in[(2N)^{-1},\infty[\cap \Q: L^{-a}\leq \eta_k\}$$
condition (\ref{lasttech1}) fails for all $\mu\in[N^{-1},N\wedge \mu_k[$, and
so (\ref{lasttech2}) must hold for the same values of $\mu$. Since (\ref{lasttech3}) shows that
$\mu_k\uparrow -\un W$ as $k\uparrow\infty$, $\N_0$ a.e. on $\{\un W <-(2N) ^{-1}\}$, this completes
the proof of Lemma \ref{occuplast} and Lemma \ref{occuptree}. \cq 

\section{Hausdorff dimension}

In this section we compute the Hausdorff dimension of the limiting metric space
appearing in Theorem \ref{main}. Although the metric $D$ is not known explicitly, 
it turns out that we have enough information to determine this
Hausdorff dimension.

\begin{theorem}
\label{Hausdim}
We have a.s.
$${\rm dim}(\t_{\ov\eg}\,/\!\approx,D)=4.$$
\end{theorem}

\proof We first derive the upper bound ${\rm dim}(\t_{\ov\eg}\,/\!\approx,D)\leq 4$. 
Recall that the process $(Z_t)_{t\in[0,1]}$ is Gaussian conditionally given $(\eg_t)_{t\geq 0}$,
and that the conditional second moment of $Z_t-Z_s$ is $m_\eg(s,t)$. Also recall that the
function $t\la \eg_t$ is a.s. H\"older continuous with exponent $\frac{1}{2}-\varepsilon$, for any 
$\varepsilon>0$. From this fact and an application of the classical Kolmogorov lemma, we
get that the mapping $t\la Z_t$ is a.s. H\"older continuous with exponent $\frac{1}{4}-\varepsilon$, for
any 
$\varepsilon\in]0,\frac{1}{4}[$. Clearly the same holds if $Z$ is replaced by $\ov Z$. Hence, if
$\varepsilon\in]0,\frac{1}{4}[$ is fixed, there exists a (random) constant $C_1$ such that, for every
$s,t\in[0,1]$,
$$|\ov Z_s-\ov Z_t|\leq C_1\,|s-t|^{\frac{1}{4}-\varepsilon}.$$
It immediately follows that, for every $s,t\in[0,1]$,
$$D^\circ(s,t)\leq 2C_1\,|s-t|^{\frac{1}{4}-\varepsilon}.$$
Since $D\leq D^\circ$, we see that the canonical projection from $[0,1]$
onto $[0,1]/\approx$ (equipped with the metric $D$) is H\"older continuous with exponent
$\frac{1}{4}-\varepsilon$. It follows that ${\rm dim}([0,1]\,/\!\approx,D)\leq
(\frac{1}{4}-\varepsilon)^{-1}$ and since $\varepsilon$ was arbitrary, we
get ${\rm dim}(\t_{\ov\eg}\,/\!\approx,D)={\rm dim}([0,1]\,/\!\approx,D)\leq 4$.

The proof of the corresponding lower bound requires the following lemma. Recall that
$\lambda$ denotes the uniform probability measure on $\t_{\ov\eg}$ (cf subsection 2.4). For every 
$a\in\t_{\ov\eg}$ and every $\varepsilon>0$, we set
$B_D(a,\varepsilon)=\{b\in \t_{\ov\eg}:D(a,b)<\varepsilon\} $.

\begin{lemma}
\label{Hauslem}
There exists a constant $C$ such that, for every $r\in]0,1]$,
$$E\Big[\int_{\t_{\ov\eg}} \lambda(da)\,\lambda(B_D(a,r))\Big]\leq C\,r^4.$$
\end{lemma}

Assume that the result of the lemma holds, and fix $\varepsilon\in]0,1]$. From the bound
of the lemma, we get that, for every integer $k\geq 1$,
$$E[\lambda(\{a\in\t_{\ov\eg} : \lambda(B_D(a,2^{-k}))\geq 2^{-k(4-\varepsilon)}\})]
\leq C\,2^{-k\varepsilon}.$$ 
By summing this estimate over $k$, we obtain
$$\limsup_{k\to\infty} \frac{\lambda(B_D(a,2^{-k}))}{2^{-k(4-\varepsilon)}}\leq 1\ ,\quad
\lambda(da)\hbox{ a.e.,\ a.s.}$$
By standard density theorems for Hausdorff measures, this implies that
${\rm dim}(\t_{\ov\eg}\,/\!\approx,D)\geq 4-\varepsilon$, a.s., which completes the
proof of Theorem \ref{Hausdim}. It only remains to prove Lemma \ref{Hauslem}. \cq

\smallskip
\noindent{\bf Proof of Lemma \ref{Hauslem}:} We rely on the case $k=2$
of Proposition \ref{finite-margi}.
With the notation of this proposition, we have
\ba
E\Big[\int_{\t_{\ov\eg}} \lambda(da)\,\lambda(B_D(a,r))\Big]
&=&E\Big[\int_{\t_{\ov\eg}\times\t_{\ov\eg}} \lambda(da)\lambda(db)\,
{\bf 1}_{\{D(a,b)<r\}}\Big]\\
&=&P[D(Y^\infty_1,Y^\infty_2)<r]\\
&\leq&\liminf_{n\to \infty} P\Big[d_n(Y^n_1,Y^n_2)<(4p(p-1)/9)^{1/4} n^{1/4}r
\Big].
\ea

On the other hand, it follows from Proposition \ref{pointedmap} that
\ba
P\Big[d_n(Y^n_1,Y^n_2)<(4p(p-1)/9)^{1/4} n^{1/4}r
\Big]\!\!\!\!&=&\!\!\!\!E\Big[\frac{1}{(p-1)n+2}\,\#B_n(Y^1_n,(4p(p-1)/9)^{1/4}n^{1/4}r)\Big]\\
&\build{\la}_{n\to\infty}^{}&E[\ov{\cal I}([0,r])].
\ea
Therefore we have obtained the bound
$$E\Big[\int_{\t_{\ov\eg}} \lambda(da)\,\lambda(B_D(a,r))\Big]\leq E[\ov{\cal
I}([0,r])].$$ 
Recall the notation of the proof of Lemma \ref{estimateISE} in Section 5. We know
that
$\ov{\cal I}([0,r])$ has the same distribution as ${\cal J}(r)$ under
$\N^{(1)}_0$. Furthermore the estimates recalled in the proof of 
Lemma \ref{estimateISE} imply that, for every $r\in]0,1]$,
$$\N_0\Big({\bf 1}_{\{\sigma>1/2\}}{\cal J}(r)\Big)\leq C'\,r^4$$
for a certain constant $C'$. A simple scaling argument then gives, with another constant $C$,
$$E[\ov{\cal
I}([0,r])]=\N^{(1)}_0({\cal J}(r))\leq C\,r^4.$$
This completes the proof of Lemma \ref{Hauslem}. \cq

\end{document}